\documentclass{aptpub}
\authornames{Alexander Stolyar} % insert the authors here for use in running head. If three or more authors please use (for example) M.~YARROW {\it et al}. Author names should follow the same M.~YARROW format and if two authors, separate by 'AND'.
\shorttitle{Large-scale particle system with mean-field interaction} % insert short title here for use in running head

\usepackage{amsmath}
\usepackage{amsfonts}
\usepackage{amssymb}

\usepackage{epsfig}

\usepackage{subfigure}
\usepackage{graphics}
\usepackage{graphicx}
\usepackage{color}
\usepackage{soul}
\usepackage[makeroom]{cancel}

\usepackage{latexsym}
\usepackage{graphics}
\usepackage{xspace}
\usepackage{color}
\usepackage{cite}
\usepackage{latexsym}
\usepackage{graphics}
\usepackage{xspace}
\usepackage{epsfig}
\usepackage{comment}

\newcommand{\beql}[1]{\begin{equation}\label{#1}}
\newcommand{\eeql}{\end{equation}}
\newcommand{\eqn}[1]{(\ref{#1})}

\newcommand{\R}{\mathbb{R}}
\newcommand{\pr}{\mathbb{P}}
\newcommand{\E}{\mathbb{E}}

\numberwithin{equation}{section}  % If you number theorems, etc. within sections,
\begin{document}%\recd{}{}%Do not alter this line.

\title{Large-scale behavior of a particle system with mean-field interaction: 
Traveling wave solutions} % insert title - use \\ if it requires more than one line.

\authorone[University of Illinois at Urbana-Champaign]{Alexander Stolyar} 
%Affiliation is just the name of your university or institution, for example 'University of Sheffield'. Author names should be of the form 'Mark Yarrow'. 
%Authors should be ordered alphabetically subject to the convention in that particular authors country. For example 'Remco van der Hofstad' would be listed under 'H' as is standard in the Netherlands. 

%Please use the following format for addresses and emails. The APT office will sort this out after you submit your files.
\addressone{ISE Department and Coordinated Science Lab, University of Illinois at Urbana-Champaign, Urbana, IL 61801} % Your postal address goes here.
\emailone{stolyar@illinois.edu} %Authors email goes here.

\begin{abstract}
We use probabilistic methods to study properties of mean-field models, arising as large-scale limits of certain particle systems with mean-field interaction. The underlying particle system is such that $n$ particles move forward on the real line.
Specifically, each particle ``jumps forward''  at some time points,
with the instantaneous rate of
jumps given by a decreasing function of the particle's location quantile within the overall distribution of particle locations.
A mean-field model describes the evolution of the particles' distribution, when $n$ is large. It
is essentially a solution to an integro-differential equation within a certain class.
Our main results concern the existence and uniqueness of -- and attraction to -- mean-field models which are traveling waves,
under general conditions on the jump-rate function and the jump-size distribution.
\end{abstract}

%insert keywords separated by a semicolon. You should avoid including keywords which also appear in the title.

% insert the primary 2020 Maths Subject Classification number in the first bracket
		% and the secondary ams number(s) in the second bracket
		% e.g. \ams{60E20}{49G03;49F10}
		%Maximum of three in each, ideally one or two in each primary and secondary.
		%codes found here ``https://mathscinet.ams.org/msnhtml/msc2020.pdf''
		
		\section{Introduction}

In this paper we use probabilistic methods to study properties of mean-field models, describing large-scale behavior of certain particle systems with mean-field interaction. A mean-field model is essentially a solution to an integro-differential equation within a certain class.
Our focus is on the existence and uniqueness of -- and attraction to -- mean-field models which are traveling waves. 

\subsection{A particle system giving rise to the mean-field model.}

The basic particle system (or, rather, its special case) which gives rise to our mean-field model
was first introduced and studied in \cite{GSS96, GMP97}, and
is as follows.
There are $n$ particles, moving in the positive direction (``right'') on the real axis $\R$. Each particle moves in jumps, as follows. For each particle there is an independent Poisson process of rate $\mu>0$. At the time points of this Poisson process the particle jumps to the right with probability $\eta_n(\nu)$, where $\nu$ is its quantile in the current empirical distribution of the particles' locations;
that is, $\nu=\ell/n$ when the particle location is $\ell$-th from the left. With complementary probability $1-\eta_n(\nu)$ the particle does not jump.
To have the model well-defined, assume that quantile-ties between co-located particles are broken uniformly at random. 
Assume that, for each $n$, $\eta_n(\nu), ~0\le \nu \le 1,$ is non-increasing, and that, as $n\to\infty$, it uniformly converges 
to a continuous, strictly decreasing function $\eta(\nu), ~0\le \nu \le 1,$ 
with $\eta(0) = 1$, $\eta(1) = 0$. The jump sizes, when a particle does jump, are given by i.i.d. non-negative r.v. with CDF $J(y), y\ge 0$; we denote by $\bar J(y) = 1- J(y)$ the complementary CDF.
In this paper we assume that for some integer $\ell \ge 2$, the jump size distribution has finite $\ell$-th moment and denote
\beql{eq-moment-finite}
m^{(k)} \doteq \int_0^\infty y^k dJ(y) < \infty, ~~k=1,2,\ldots,\ell.
\eeql
So, $m^{(1)}<\infty$ is the mean jump size.

Let $f^n(x,t)$ be the (random) empirical distribution of the 
particle locations at time $t$; namely, $f^n(x,t)$ is the fraction of particles located in $(-\infty,x]$ at time $t$.
As $n\to\infty$, it is very intuitive that $f^n(x,t)$ converges (in appropriate sense, under appropriate conditions) to a deterministic 
function $f(x,t)$, such that $f(\cdot,t)$ is a distribution function for each $t$, and the following equation holds:
\beql{eq-dyn-trans-intro}
\frac{\partial}{\partial t} f(x,t) =
- \mu \int_{-\infty}^x d_y f(y,t) \eta(f(y,t)) \bar J(x-y),
\eeql
where $d_y$ means the differential in $y$. 

At this point let us describe only the intuition for \eqn{eq-dyn-trans-intro}. 
For each $t$, the distribution $f(\cdot,t)$ approximates the distribution of particles $f^n(\cdot,t)$ when $n$ is large.
Since particles move right, $f(x,t)$ is non-increasing in $t$ for each $x$. So, $(\partial/\partial t) f(x,t) \le 0$ and its value should be equal to
the RHS of \eqn{eq-dyn-trans-intro}, which gives the instantaneous rate (scaled by $1/n$ and taken with minus sign) at which particles jump over point $x$ at time $t$. We will call a function $f(x,t)$ satisfying \eqn{eq-dyn-trans-intro} a {\em mean-field model.} 
The formal meaning of \eqn{eq-dyn-trans-intro} and the definition of a mean-field model will be given later.

It is also intuitively clear (and easy to make formal, as we do later) that, for any mean-field model,
 the speed, at which the mean $\int x d_x f(x,t)$
of the distribution $f(\cdot,t)$ moves right, must be equal to $v = m^{(1)} \mu  \int_0^1 \eta(\nu) d\nu.$ 
Suppose a mean-field model $f(x,t)$ that is a traveling wave exists, namely $f(x,t) = \phi(x-vt)$ 
for some distribution function $\phi(\cdot)$ which we will call a {\em traveling wave shape}.
By substituting into \eqn{eq-dyn-trans-intro}, we see that any traveling wave shape $\phi$
must satisfy equation
\beql{eq-wave}
v  \phi'(x) =
\mu \int_{-\infty}^x \phi'(y) \eta(\phi(y)) \bar J(x-y) dy.
\eeql
We will make this statement formal later. 

\subsection{Motivation for the particle system.}
\label{sec-motivation-particle}

The original motivation for the described particle system is an idealized model \cite{GSS96, GMP97}
of distributed parallel simulation. The $n$ particles represent $n$ components, or sites, of one large system.
Each site is being simulated by a separate computer (processor). 
The sites are interdependent, so their simulations cannot run independently. 
A particle location at a given real time represents the ``local simulation time'' of the corresponding site, which is
the time in the simulated system up to which the simulation of this site is valid.
After a site's independent simulation runs for some real time, the site tries 
to update its state and advance its local time. However, this local time advance is not always possible, because
the site's evolution depends on a number of other sites, whose local times may be lagging behind; if the 
site's local time cannot be advanced, it will ``go back'' and start simulation again from the current local time. 
In general, a site's local time advance is more likely to occur if it is ``further behind'' the local times of other sites. 
The model in \cite{GSS96, GMP97} assumes specifically that: each 
particle gets ``urges'' to jump forward (advance local time) as an independent Poisson process of rate $\mu$;
when a particle does get a jump urge, it actually jumps only if $K$ other particles, $K \ge 1$, 
chosen uniformly at random, are currently ahead of it.
If $\nu$ is the quantile of the particle location, then the probability of the particle jump, $\eta_n(\nu)$, is such that 
$\eta_n(\nu) \to \eta(\nu) = (1-\nu)^K$ as $n\to\infty$. It is also assumed in \cite{GSS96, GMP97} that the jump size distribution $J(\cdot)$ is exponential. 
\cite{GMP97} further assumes that the system initial state is such that the locations of $n$ particles are drawn independently
from a given distribution.
We will comment on these additional assumptions in detail later, in Section~\ref{sec-prev-work}.

\subsection{Prior results on the particle system and its mean-field model.}
\label{sec-prior-brief}

Papers \cite{GSS96,GMP97} address two different issues related to the particle system and solutions of \eqn{eq-dyn-trans-intro}:
\begin{itemize}
\item {\em Limit transition from $f^n(x,t)$ to a mean-field model $f(x,t)$, as $n\to\infty$.} \\
In \cite{GMP97} it is proved, under certain additional assumptions, that, as $n\to\infty$, the random process $f^n(x,t)$ indeed converges to a deterministic process $f(x,t)$, satisfying \eqn{eq-dyn-trans-intro}.

\item {\em Convergence of a mean-field model $f(x,t)$ to a traveling wave solution, as $t\to\infty$.}\\
In \cite{GSS96} the following is proved, under certain additional assumptions. 
{\em If a traveling wave shape $\phi$, i.e. a solution to \eqn{eq-wave}, exists}, then
it is unique (up to a shift) and,  
as $t\to\infty$, 
a mean-field model $f(x,t)$ converges to the traveling wave solution, namely
$f(\cdot+vt,t) \to \phi(\cdot)$. The question of the existence of a traveling wave shape $\phi$ 
was left open, except 
for the case of exponential distribution $J(\cdot)$ and $\eta(\nu)=(1-\nu)^K$, when \eqn{eq-wave} is easily explicitly solvable.
\end{itemize}

A detailed discussion of the results of \cite{GSS96,GMP97} is given in Section~\ref{sec-prev-work}, which also discusses other related work.

\subsection{Main results of this paper.}

We study the properties of mean-field models, specifically the existence and uniqueness of, and convergence to, traveling wave shapes.
The convergence to a traveling wave shape is important, because if it holds, it means that the particles' locations ``remain close
to each other'' (as they move right) {\em regardless of the number of particles,} in the sense that distribution of the particles' locations stays close to a certain shape, which moves right at the constant speed $v$.

Our main results are:
\begin{itemize}
\item We prove (in Theorem~\ref{thm-existence}) the existence of a traveling wave shape $\phi$, 
for general jump size distribution $J(\cdot)$ and general (strictly decreasing continuous) $\eta(\cdot)$.
Moreover, as a distribution, $\phi$ has finite $(\ell-1)$-th moment, $\int_y |y|^{\ell-1} d\phi(y) < \infty$, for the $\ell\ge 2$ in \eqn{eq-moment-finite}.
\item Under the additional condition that $J(\cdot)$ has positive density (bounded away from $0$ on compact sets),
we show (in Theorem~\ref{thm-main}) the uniqueness (up to a shift) of the traveling wave shape $\phi$ and convergence to it, $f(\cdot+vt,t) \to \phi(\cdot)$,
for any mean-field model.
\item As the main analysis tool, we introduce and study the properties of traveling wave shapes within ``finite frames.'' 
The existence of a traveling wave shape is then proved by letting the frame size go to infinity.
These results may be of independent interest. 
\end{itemize}

We emphasize that the results of this paper concern properties of formally defined mean-field models. The question of whether/when the convergence to a mean-field model holds, is a separate issue, partially addressed in \cite{GMP97}, 
as pointed out above in Section~\ref{sec-prior-brief} and will be discussed in more detail in Section~\ref{sec-prev-work-conv-to-mf}.

\subsection{Overview of the technical approach and key challenges.}

Our approach to proving a traveling wave shape existence (Theorem~\ref{thm-existence}) is as follows. A traveling wave shape $\phi$
can be characterized as a fixed point of an operator $A$, which maps a probability distribution $\phi$ (describing a distribution of particle locations) into the stationary distribution $A \phi$ of a single particle,
evolving within the ``environment'' given by $\phi$. 
(This is done in Section~\ref{sec-wave-charcter}.)
Specifically, the particle jumps right as described in our motivating particle system, except its current location quantile 
is that within distribution $\phi$; between jumps the particle moves left at the constant speed $v$. 
Since $\phi$ is a fixed point, it is natural to try to obtain its existence using in some way the Brouwer fixed point theorem. 
In our case this cannot be done directly, because the space of proper distributions on the real line in non-compact. 

To address this challenge, we introduce and study traveling wave shapes within ``finite frames'' --
they are fixed points of finite-frame versions of the operator $A$, and play a key role in our analysis.
(This is done in Section~\ref{sec-wave-finite-frame}.)
A finite-frame version $A^{(w;B)}$ of the operator $A$ maps a distribution $\gamma$, concentrated on the compact 
segment (finite frame) $[-B,B]$,
into the stationary distribution $A^{(w;B)} \gamma$ of a single particle,
which evolves within the environment given by $\gamma$, and whose movement is restricted to $[-B,B]$ as follows: 
the particle cannot "jump over $B$;'' between jumps it moves left at the constant speed $w$, but cannot move to the left of $-B$.
We show that operator $A^{(w;B)}\gamma$ is continuous in $\gamma$, as well as in the parameters $B$ and $w$.
Given the compactness of a finite frame, the existence of fixed points of  $A^{(w;B)}$  -- finite-frame traveling wave shapes -- follows.
We further show that a finite-frame traveling wave shape is unique and depends continuously on parameters $B$ and $w$;
showing the uniqueness (Lemma~\ref{lem-fp-new}(i)) is the most involved part here. 

We then let $B\to\infty$ along some sequence, and for each $B$ we choose $w=w(B)$ such that the corresponding finite-frame traveling wave shape $\gamma_B$ has its median exactly at $0$. The rest of a traveling wave shape existence proof 
(in Section~\ref{sec-main-proofs-existence}) shows that the family of distributions $\{\gamma_B\}$ is tight, and any subsequential limit $\phi$
must be a traveling wave shape. Proving the tightness is the key technical part of the paper. It involves, first, 
showing (Lemma~\ref{lem111}) that, necessarily, $w\to v$ as $B\to\infty$; if not, the particle process would have inherent positive or negative drift, and the proof obtains a contradiction with the fact that the median of $\gamma_B$ stays at $0$. Finally, Lemma~\ref{lem-gamma_B-tight} 
shows the tightness of $\{\gamma_B\}$ itself; if not, the proof considers the single-particle process under a space as well as a space/time rescaling and shows that, for large $B$, these rescaled processes would have a strictly negative steady-state drift of a quadratic Lyapunov function -- a contradiction. The use of quadratic Lyapunov function here requires that the jump size distribution has finite second moment;
this is where the finite second moment assumption is used in the existence proof.

\subsection{Outline of the rest of the paper.}

Section~\ref{sec-notation} gives some basic notation used throughout the paper.
In Section~\ref{sec-main-results} we formally define mean-field models, state of main results, and 
discuss previous related work. Results in Section~\ref{sec-wave-charcter} characterize a traveling wave shape 
as a fixed point an operator $A$; here we also obtain the finiteness of moments of $A \phi$.
In Section~\ref{sec-wave-finite-frame} we introduce and study traveling wave shapes within finite frames.
Section~\ref{sec-main-proofs} contains the proofs of our main results, Theorems~\ref{thm-existence} and \ref{thm-main}. The discussion in Section~\ref{sec-discusion} includes: a generalization of our main results; possible relaxation of assumptions; and
a conjecture about the limit of the stationary distributions of the particle system, as the number of particles $n \to \infty$.

\section{Basic notation}
\label{sec-notation}

Sets of real and real non-negative numbers are denoted by $\R$ and $\R_+$, respectively. 
As a measurable space, $\R$ is endowed with Borel $\sigma$-algebra $\mathcal B(\R)$. 
Convergence to a set, $x\to S \subseteq \R$,
means $\inf_{y\in S} |x-y| \to 0.$ For $x\in \R$, $\lfloor x \rfloor$ is the largest integer not greater than $x$.
For a condition/event $A$, $I\{A\}=1$ if $A$ holds, and $I\{A\}=0$ otherwise.

For functions $h(x)$ of a real $x$: $h(x+)$ and $h(x-)$ are the right and left limits;
$(d/dx) h(x)=h'(x)$ is the derivative; $(d^+ / dx) h(x)$ and 
$(d^- / dx) h(x)$ are the right and left derivatives;
$$
\frac{d^+_\ell}{dx} h(x) = \liminf_{y\downarrow x} \frac{h(y)-h(x)}{y-x}, ~~\frac{d^+_u}{dx} h(x) = \limsup_{y\downarrow x} \frac{h(y)-h(x)}{y-x}, 
$$
are the lower and upper right derivative numbers; $\|h\| = \sup_x |h(x)|$ is the sup-norm, and $\stackrel{u}{\rightarrow}$ is the corresponding uniform convergence; $\|h\|_1 = \int_x |h(x)|dx$ is $L_1$-norm and $\stackrel{L_1}{\rightarrow}$ is the corresponding convergence;  $h \stackrel{u.o.c.}{\rightarrow} g$ means {\em uniform on compact sets} (u.o.c.) convergence; we denote by $\theta_c, ~c\in\R$, the shift operator $\theta_c h(x) = h(x-c)$; $h(x)$ is called $c$-Lipschitz if it is Lipschitz with constant $c \ge 0$. Notation $d_x h(x,t)$ for a multivariate function
 means the differential in variable $x$.

We say that a function $g(x)$ of a real $x$ is RCLL if it is right-continuous with left limits. The domain of $g(x)$ will be clear from the context; usually, $x\in \R$. For RCLL functions, $\stackrel{J_1}{\rightarrow}$ denotes Skorohod ($J_1$) convergence (cf. \cite{Ethier_Kurtz}).

For non-decreasing RCLL functions, $h \stackrel{w}{\rightarrow} \gamma$ denotes the weak convergence, namely, the
convergence at every point of continuity of $\gamma$. Symbol $\stackrel{w}{\rightarrow}$ is also used more generally, to denote 
weak convergence of measures.

A non-decreasing RCLL function $\gamma=(\gamma(x), ~x\in \R)$ is a probability distribution function, if $\lim_{x\downarrow -\infty} \gamma(x) \ge 0$ and $\lim_{x\uparrow \infty} \gamma(x) \le 1$; a probability distribution function $\gamma$ is {\em proper} if $\lim_{x\downarrow -\infty} \gamma(x) = 0$ and $\lim_{x\uparrow \infty} \gamma(x) = 1$; thus, an improper $\gamma$ may have atoms at $-\infty$ and/or $\infty$.
We use the terms {\em probability distribution function}, {\em distribution function} and {\em distribution} interchangeably. 
Unless explicitly stated otherwise, a distribution means proper distribution. The inverse ($\nu$-th quantile) of a 
(proper or improper) distribution $\gamma$ is $\gamma^{-1}(\nu) \doteq \inf\{y~|~\gamma(y) \ge \nu\}$, $\nu\in [0,1]$; $\gamma^{-1}(\nu) =\infty$ when the set under $\inf$ is empty. 
We use a usual stochastic order (dominance) relation between probability distributions on $\R$: 
$g \preceq \gamma$  iff $g(x) \ge \gamma(x), ~x \in \R$; we refer to this as $\gamma$ dominating 
%(or being larger than, or being ahead of, or being to the right of) 
$g$. 
For a distribution $\gamma$ and a function $h$: $\gamma h \doteq \int_\R h(y) d\gamma(y)$,
$$
\int_{x-}^\infty h(y) d\gamma(y) \doteq \int_{[x,\infty]} h(y) d\gamma(y), ~~
\int_{x+}^\infty h(y) d\gamma(y) \doteq \int_{(x,\infty]} h(y) d\gamma(y).
$$

When $G_k, G$ are operators mapping a function of $x\in \R$ into another function of $x\in \R$, the convergence $G_k h\to Gh$, 
or $\lim G_k h = Gh$,
always means uniform convergence $\| G_k h - G h \| \to 0$.

Suppose we have a Markov process taking values in $\R$, with $P^t(x,H),$ $t\ge 0, ~x\in \R, ~H\in \mathcal B(\R),$ being its transition function. $P^t$, as an operator, is $P^t h(x)  \doteq \int_y P^t(x,dy) h(y)$; $I=P^0$ is the identity operator.
The process (infinitesimal) generator $G$ is 
$$
Gh \doteq \lim_{t\downarrow 0} (1/t) [P^t - I] h.
$$
Function $h$ is within the domain of the generator $G$ if $Gh$ is well-defined.

We will also use the following non-standard notation throughout the paper. 
For a probability distribution function $\gamma$, a strictly decreasing continuous function $\eta(\nu), ~0\le \nu \le 1$,   
and $y \in \R$, we denote 
$$
    \bar \eta(y,\gamma) \doteq \left\{\begin{array}{ll}
        \eta(\gamma(y)), & \text{when $\gamma(\cdot)$ is continuous at point $y$}\\
        (\nu_2-\nu_1)^{-1} \int_{\nu_1}^{\nu_2} \eta(\nu) d\nu, & \text{otherwise},
        \end{array}\right. 
$$
where $\nu_2=\gamma(y), \nu_1=\gamma(y-)$.

Convergence in probability is denoted by $\stackrel{P}{\longrightarrow}$;
{\em w.p.1} means {\em with probability $1$};  
RHS and LHS mean right-hand side and left-hand side, respectively; WLOG means {\em without loss of generality}.

\section{Mean-field model of a large-scale system}
\label{sec-main-results}

\subsection{Mean-field model definition}

Let $\eta(\nu), ~0\le \nu \le 1,$ be a a continuous, strictly decreasing function,
with $\eta(0) = 1$, $\eta(1) = 0$. Let $J(y), y\ge 0$, be the CDF of a probability distribution, concentrated on $\R_+$ and satisfying \eqn{eq-moment-finite}
for some integer $\ell \ge 2$; denote $\bar J(y) = 1- J(y)$.
We now introduce the following

\begin{defn}
\label{def-mfm}
A function $f(x,t), ~x\in \R, ~t\in \R_+,$ will be called a {\em mean-field model} if it satisfies the following conditions.\\
 (a) For any $t$, as a function of $x$, $f(x,t)$ is a probability distribution function; that is $f(\cdot,t)$ is non-decreasing RCLL, with 
$\lim_{x\to -\infty}f(x,t) =0$ and $\lim_{x\to \infty}f(x,t) =1$.\\
(b) For any $x$, $f(x,t)$ is non-increasing $c$-Lipschitz in $t$, with constant $c$ independent of $x$.\\
(c) For any $x$, for any $t$ where the partial derivative
$(\partial/\partial t) f(x,t)$ exists (which is almost all $t$ w.r.t. Lebesgue measure, by the Lipschitz property), equation 
\beql{eq-dyn-trans}
\frac{\partial}{\partial t} f(x,t) =
- \mu \int_{-\infty}^x d_y f(y,t) \bar \eta(y,f(\cdot,t)) \bar J(x-y),
\eeql
holds.
\end{defn}

Equation \eqn{eq-dyn-trans} is a more general form of \eqn{eq-dyn-trans-intro}, allowing $f(x,t)$ to be RCLL in $x$, rather than continuous.
If $f(\cdot,t)$ is continuous at $y$, then  $\bar \eta(y,f(\cdot,t))=\eta(f(y,t))$;
if $f(\cdot,t)$ has a jump at $y$ then $\bar \eta(y,f(\cdot,t))$ is
$\eta(\nu)$ averaged over $\nu\in [f(y-,t),f(y,t)]$.

Note that the following holds for any mean-field model $f(x,t)$. Denote
$$
v \doteq \mu m^{(1)} \int_0^1 \eta(\nu) d\nu.
$$
Then for any $\tau \le t$,
\beql{eq-speed-conserv} 
\int_{-\infty}^{\infty} [f(x,\tau) - f(x,t)] dx = v (t-\tau).
\eeql
Indeed, if we denote by $h(x,t)$ the RHS of \eqn{eq-dyn-trans}, we obtain
$$
\|h(\cdot,t)\|_1 = 
 \int_{-\infty}^{\infty} dx ~\mu \int_{-\infty}^x d_y f(y,t) \bar \eta(y,f(\cdot,t)) \bar J(x-y) =
$$
$$
 \int_{-\infty}^{\infty} dx ~\mu \int_{0}^{f(x,t)} d\nu ~\eta(\nu)   \bar J(x-f^{-1}(\nu,t)) = 
$$
$$
 \mu \int_{0}^{1} \eta(\nu) d\nu 
 \int_{f^{-1}(\nu,t)}^{\infty} dx ~  \bar J(x-f^{-1}(\nu,t)) = 
$$
$$
 \mu \int_{0}^{1} \eta(\nu) d\nu 
 \int_{0}^{\infty} \bar J(\xi) d\xi = \mu m^{(1)} \int_0^1 \eta(\nu) d\nu =v,
$$
where $f^{-1}(\nu,t)$ is the inverse of $f(y,t)$ with respect to $y$. For any $x$, since $f(x,t)$ is Lipschitz -- and then absolutely continuous -- in $t$, we have
$$
f(x,\tau) - f(x,t) = \int_{\tau}^t [-h(x,\xi)] d\xi.
$$
Then
$$
\int_{-\infty}^{\infty} [f(x,\tau) - f(x,t)] dx = \int_{-\infty}^{\infty} \int_{\tau}^t [-h(x,\xi)] d\xi dx =
$$
$$
 \int_{\tau}^t \int_{-\infty}^{\infty} [-h(x,\xi)]  dx d\xi = \int_{\tau}^t \| h(\cdot,\xi)\|_1 d\xi = v (t-\tau).
$$
Equality (``conservation law'') \eqn{eq-speed-conserv} implies in particular that, if the mean of the distribution $f(\cdot,\tau)$,
$$
\bar f(\tau) \doteq \int_{-\infty}^{\infty} x d_x f(x,\tau),
$$
is well-defined and finite (i.e., $\int_{-\infty}^{\infty} |x| d_x f(x,\tau) < \infty$), then $\bar f(t)$ is finite for all $t\ge\tau$, and
$$
\bar f(t) - \bar f(\tau) = v(t-\tau).
$$
In other words, if $\bar f(\tau)$ is finite, then $\bar f(t)$ is finite for all $t\ge\tau$, 
it moves right at the constant speed $v$.

\begin{defn}
\label{def-tws}
Suppose a mean-field model $f(x,t)$, which is a traveling wave, exists; namely $f(x,t) = \phi(x-vt)$ for some probability distribution function
$\phi(\cdot)$. Such a function $\phi$ will be called a {\em traveling wave shape}.
\end{defn}

Note that by \eqn{eq-speed-conserv} the speed of a traveling wave can only be $v$.
It is straightforward to see that {\em a function $\phi$ is a traveling wave shape if and only if it is 
a Lipschitz continuous distribution function,
satisfying \eqn{eq-wave} at any $x$ where $\phi'(x)$ exists (which is almost any $x$).}
Furthermore, from the form of \eqn{eq-wave} we see that, in fact,
{\em $\phi$ satisfies \eqn{eq-wave} for each $x$,  and the derivative $\phi'(x)$ is continuous.}
\iffalse
Substituting $f(x,t) = \phi(x-vt)$ into \eqn{eq-dyn-trans}, we observe that any traveling wave shape $\phi(x)$ must be Lipschitz continuous,
and then \eqn{eq-dyn-trans} in fact takes a simpler form \eqn{eq-dyn-trans-intro}. 
Therefore, a traveling wave shape $\phi$ is a Lipschitz continuous distribution function,
satisfying \eqn{eq-wave} at almost any $x$. Then, from the form of \eqn{eq-wave} we see that, in fact,
$\phi$ satisfies \eqn{eq-wave} for each $x$,  and the derivative $\phi'(x)$ is continuous.
Conversely, any Lipschitz continuous distribution function $\phi$,
satisfying \eqn{eq-wave} at almost any $x$ (and then necessarily at any $x$), is a traveling wave shape. 
\fi

Note that if $\phi(x)$ is a traveling wave shape, then so is $\phi(x-c)$ for any constant real shift $c$. 

\subsection{Main results.}

The following is the main result of this paper. It proves the existence of a traveling wave shape $\phi$, under the general assumptions of our model. It also shows that $\phi$ has a finite $(\ell-1)$-th moment.

\begin{thm}
\label{thm-existence} Assume \eqn{eq-moment-finite}. Then:

(i) There exists a traveling wave shape $\phi(\cdot)$. 

(ii) Any traveling wave shape $\phi(\cdot)$ is such that, for the integer $\ell \ge 2$ in the condition \eqn{eq-moment-finite},
\beql{eq-moments}
\int_0^\infty |y|^{\ell-1} d\phi(y) < \infty.
\eeql
\end{thm}

We also obtain the following uniqueness and convergence result. It is proved by
combining Theorem~\ref{thm-existence} with the results of \cite{GSS96}.

\begin{thm}
\label{thm-main}
Assume, in addition, that the jump size distribution is such that
\beql{eq-positive-density}
\mbox{$J(\cdot)$ has density $J'(y)>0$, bounded away from $0$ on compact subsets of $\R_+$}.
\eeql
Then the following holds.

(i) The traveling wave shape $\phi(\cdot)$ is unique, up to a shift $\phi(\cdot-c)$.

(ii) If  a mean-field model $f(x,t)$ is such that the initial mean
$\bar f(0)$ is well-defined and finite, then the convergence to the unique traveling wave shape $\phi$ takes place,
\beql{eq-conv-to-shape}
f(\cdot+vt,t) \stackrel{L_1}{\rightarrow} \phi(\cdot)  ~~\mbox{and}~~   f(\cdot+vt,t) \stackrel{u}{\rightarrow} \phi(\cdot),
\eeql
where $\phi$ is uniquely centered by condition $\int_y y d\phi(y)=\bar f(0)$. 
\end{thm}

Note that to have a ``clean'' convergence to the traveling wave shape $\phi$, as in \eqn{eq-conv-to-shape},
{\em some} additional conditions on the distribution $J(\cdot)$ 
are required. For example, if $f(\cdot,0)$ is concentrated on a lattice $\{c k, ~\mbox{$k$ is integer}\}$,
and distribution $J(\cdot)$ is arithmetic, concentrated on $\{c k, ~k=1,2,\ldots\}$, 
then the convergence \eqn{eq-conv-to-shape} is impossible, even though a traveling wave shape $\phi$ does exist
(by Theorem~\ref{thm-existence}), and might be unique.

Most of the rest of the paper is devoted to the proof of Theorem~\ref{thm-existence}. The proof constructs a traveling wave shape $\phi$ as a limit of traveling wave shapes within finite frames. The finite-frame traveling wave shapes and their analysis might be of independent interest.

Before proceeding with the proofs, we observe that, WLOG, we can assume that
\beql{eq-mu-m-wlog}
\mu=1, ~~m^{(1)}=1.
\eeql
Indeed, any mean-field model can be reduced to the corresponding model satisfying \eqn{eq-mu-m-wlog}, by time and space rescaling.
So, in all the proofs we do assume \eqn{eq-mu-m-wlog},
in which case 
$$
v=\int_0^1 \eta(\xi) d\xi,
$$
\eqn{eq-dyn-trans} becomes 
\beql{eq-dyn-trans2}
\frac{\partial}{\partial t} f(x,t) =
- \int_{-\infty}^x d_y f(y,t) \eta(f(y,t)) \bar J(x-y), 
\eeql
and \eqn{eq-wave} becomes
\beql{eq-wave2}
v  \phi'(x) =
\int_{-\infty}^x \phi'(y) \eta(\phi(y)) \bar J(x-y) dy.
\eeql

\subsection{Previous work}
\label{sec-prev-work}

\subsubsection{Detailed discussion of \cite{GMP97}.}
\label{sec-prev-work-conv-to-mf}

As mentioned earlier, paper \cite{GMP97} proves that, as $n\to\infty$, the random process $f^n(x,t)$ converges to a mean-field model $f(x,t)$,
{\em under the additional assumptions} that we now discuss. 

The first additional assumption is on {\em the form of  the jump probability 
functions $\eta_n(\cdot)$ in the pre-limit particle system.}
The paper considers the specific particle system, in which when a particle gets a jump urge, it actually jumps only if $K$ other particles, 
$K\ge 1$, chosen uniformly at random, are ahead of it. In this case $\eta_n(\nu) \to \eta(\nu) = (1-\nu)^K$ 
(see Section~\ref{sec-motivation-particle}). The important part of this assumption about $\eta_n(\nu)$ is that 
the decision on whether or not the particle jumps depends only on the locations of a uniformly at random chosen set of other particles,
and the cardinality of that random set is upper bounded by some $K$. Consequently, the jump probability assumption in \cite{GMP97}
can be generalized, for example, as follows: a particle that gets an urge to jump chooses $K$ other particles uniformly at random and its decision to jump or not depends only on $k/K$, where $k$ is the number of chosen particles that are behind this particle. 
Under this more general assumption on $\eta_n(\nu)$, the results of \cite{GMP97} still apply, as long as $\eta_n(\cdot)$
uniformly converges to some strictly decreasing continuous function $\eta(\cdot)$, with $\eta(0) =1$
and $\eta(1) =0$. Now, suppose we a priori fix some arbitrary strictly decreasing continuous function $\eta(\nu)$, with $\eta(0) =1$
and $\eta(1) =0$. Fix integer $K\ge 1$. Consider the following jump probability rule (which uniquely determines the corresponding
jump probability function $\eta^{(K)}_n(\cdot)$): a particle which gets an urge to jump chooses $K$ other particles uniformly at random and
 actually jumps with (random) probability $\eta(k/K)$, where $k\ge 0$ is the (random) number of chosen particles that happen to be behind our particle. It is easy to see that, as $n\to\infty$,
$$
\eta^{(K)}_n(\nu)  \to \eta^{(K)}(\nu) \doteq \sum_{k=0}^K \frac{K!}{k!(K-k)!} \nu^k (1-\nu)^{K-k} \eta(k/K).
$$
The results of \cite{GMP97} apply to such system, with its mean-field limit $f(x,t)$ having $\eta^{(K)}(\cdot)$ 
in place of $\eta(\cdot)$. Note that, if $K$ is large (but fixed!), $\eta^{(K)}(\cdot)$ can arbitrarily closely approximate $\eta(\cdot)$.
To summarize, the results of \cite{GMP97} do {\em not} apply to show the convergence to a mean-field model 
with {\em any} a priori fixed jump probability function $\eta(\cdot)$. 
(Recall that in this paper we define a mean-field model formally, for {\em any} such  function $\eta(\nu)$.) 
However, there always exists a function $\eta^{(K)}(\cdot)$, arbitrarily close to $\eta(\cdot)$,
such that the results of \cite{GMP97} apply to obtain the convergence to a mean-field model 
with $\eta^{(K)}(\cdot)$ in place of $\eta(\cdot)$. 

The second additional assumption in \cite{GMP97} is that, for any $n$, the system initial state is such that 
the locations of the $n$ particles are drawn {\em independently} from a given {\em absolutely continuous} distribution $f(\cdot,0)$
with finite second moment.

Finally, \cite{GMP97} assumes that the jump size distribution $J(\cdot)$ is exponential. The results
of \cite{GMP97} easily generalize to an arbitrary absolutely continuous distribution $J(\cdot)$.

\subsubsection{Detailed discussion of \cite{GSS96}.}

Paper \cite{GSS96} considers a formally defined mean-field model (in the terminology of this paper), under the additional assumptions 
that $J(\cdot)$ is exponential, $J(y) = 1-e^{-y}$, 
and $\eta(\nu) = (1-\nu)^K$. 

The main result of \cite{GSS96} is basically as follows. Consider two mean-field models, $f_1(x,t)$ and $f_2(x,t)$, 
with well-defined, equal means, $\bar f_1(t)=\bar f_2(t)=\bar f_1(0)+vt=\bar f_2(0)+vt$. Then, the derivative of the $L_1$-distance is negative,
$(d/dt) \|f_1(\cdot,t) - f_2(\cdot,t)\|_1 < 0$, and bounded away from $0$ as long as $\|f_1(\cdot,t) - f_2(\cdot,t)\|_1$ is bounded
away from $0$. This, in turn, leads to the second result: {\em if a traveling shape $\phi$ exists}, then it is unique (up to a shift) and 
the convergence \eqn{eq-conv-to-shape} to the traveling wave shape holds. These results, along with their proofs, extend to the 
more general case of arbitrary continuous strictly decreasing $\eta(\nu)$ and distribution
$J(\cdot)$ satisfying \eqn{eq-positive-density}. 

As far as the existence of a traveling wave shape is concerned, in \cite{GSS96} it is demonstrated {\em only
for the special case of exponential $J(\cdot)$
and $\eta(\nu) = (1-\nu)^K$}, in which case \eqn{eq-wave2} can be explicitly solved as follows. 
The derivation is not given in \cite{GSS96}, so, for completeness, we give it here.

In the case $J(y)=1-e^{-y}$, \eqn{eq-wave2} becomes an ODE. Indeed, we have
$$
v  \phi'(x) = \int_{-\infty}^x \phi'(y) \eta(\phi(y)) e^{-(x-y)} dy,
$$
$$
v  \phi'(x) e^x = \int_{-\infty}^x \phi'(y) \eta(\phi(y)) e^y dy,
$$
\beql{eq-wave2-exp}
v  [\phi''(x) + \phi'(x)] = \phi'(x) \eta(\phi(x)).
\eeql
Under the further assumption that $\eta(\nu) = (1-\nu)^K$, \eqn{eq-wave2-exp} is easy to solve explicitly. 
In this case $v=1/(K+1)$.
Denote $\alpha(x) =1-\phi(x)$.
Then
$$
v  [-\alpha''(x) - \alpha'(x)] = -\alpha'(x) \eta(1-\alpha(x)) = -\alpha'(x) (\alpha(x))^K = -\frac{[(\alpha(x))^{K+1}]'}{K+1},
$$
$$
-\alpha''(x) - \alpha'(x) = -[(\alpha(x))^{K+1}]',
$$
$$
-\alpha'(x) - \alpha(x) = -(\alpha(x))^{K+1} +c_1.
$$
Recalling that $1-\alpha(x)$ is a probability distribution function, we see that $c_1=0$ must hold, so
$$
-\alpha'(x) - \alpha(x) = -(\alpha(x))^{K+1},
$$
$$
-\frac{(\log \alpha(x))'}{1-(\alpha(x))^K} = 1,
$$
$$
\frac{1}{K} \left[ \log[(\alpha(x))^{-K} -1]\right]' = 1.
$$
Solving for $\alpha(x)$, we finally obtain
\beql{eq-wave2-exp-final}
\phi(x) = 1- \left[ 1+ e^{K(x-c)}\right]^{-1/K},
\eeql
where real $c$ is a (shift) parameter.

\subsubsection{Other related work.}

Subsequently to \cite{GMP97,GSS96}, there has been a line of work generally focused on modeling mean-field interaction 
between elements (particles) of a large system, 
cf. \cite{Man-Sch-2005, Mal-Man-2006, Manita-2009, Manita-2006, Malyshkin-2006,Manita-2014,balazs-racz-toth-2014}. The models and results in this line of work differ from this paper in that they consider
different particle jump rules and/or different asymptotic regimes. But the general common feature is that they study systems with
mean-field (``global'') interactions, when each particle's instantaneous behavior 
depends on the current distribution of all particles' states, as opposed to ``local'' interactions,'' when a particle's instantaneous behavior 
depends on the current states of its ``neighbors.'' This naturally leads to the corresponding mean-field models, describing 
behavior of large-scale systems. 

In particular, paper \cite{balazs-racz-toth-2014} considers a related particle system with the jump probability depending not on a particle quantile within the empirical distribution $f^n(\cdot,t)$, but rather on the particle displacement with respect to the empirical distribution mean. For this system, under some technical assumptions, the paper proves the convergence of $f^n(x,t)$ to the corresponding (different from ours) mean-field model $f(x,t)$. This parallels the main result of \cite{GMP97} for our particle system.
The paper also shows the existence of traveling wave shapes for that mean-field model, but only for the special case when the jump size distribution is exponential.

Further references and discussion of related work can be found in e.g. \cite{Manita-2014,balazs-racz-toth-2014}.

\section{Characterization of a traveling wave shape}
\label{sec-wave-charcter}

Consider the following operator $A$, mapping a proper probability distribution $\psi(\cdot)$ on $\R$ into another distribution $A \psi(\cdot)$. 
Consider one particle, whose location $X(t)$ evolves as follows. 
The particle moves left at the constant speed $v$ and sometimes jumps right, with i.i.d. jump sizes given by the distribution $J(\cdot)$.
The particle get urges to jump as a rate $1$ Poisson process, and actually jumps with probability 
$\bar \eta(X(t),\psi)$, depending on the current particle location $X(t)$ and distribution $\psi$.
(Since the set of discontinuity points of $\psi$ is at most countable, w.p.1 the process $X(\cdot)$ 
will never jump from points of discontinuity of $\psi$. So, equivalently, we can say that the jump probability is $\eta(\psi(X(t)))$.)
The distribution $A \psi$ is defined as the stationary distribution of the process $X(\cdot)$. The following lemma 
shows that this stationary distribution exists and is unique.

\begin{lem}
\label{lem-stability}
Consider the single-particle movement process $X(\cdot)$ defined just above. This process is positive recurrent, 
and therefore has unique stationary distribution. Moreover, 
in steady-state
\beql{eq-moments-oper}
\E |X(t)|^{\ell-1} < \infty,
\eeql
for the $\ell\ge 2$ from condition \eqn{eq-moment-finite}.
\end{lem}

\begin{proof} 
The distribution $\psi$ is proper. 
Then, when $X(t) < 0$ with large absolute value, the process drift is positive, close to 
$1-v$; when $X(t)$ is large positive, the process drift is negative, close to $-v$. The formal positive recurrence proof is a straightforward application
of the fluid limit technique \cite{RS92, Dai95, St95, Bramson-book}. Indeed, 
it is easy to see that any compact subset of $\R$ is petite; see the definition in e.g. \cite{Dai95,Bramson-book}.
(For a compact set, consider an interval $[-C,C]$ containing it. Then the sampling distribution can be chosen as, for example,
the uniform distribution in $[0,3C/v]$. Here we use the fact that in the interval $[0,3C/v]$, with probability $1-e^{-3C/v}$,
the particle does not have an urge to jump and will simply move left at the constant speed $v$, covering the distance $3C$
and in particular the entire interval $[-2C,-C]$.)
Further, if we consider a sequence of processes $X^{(r)}(\cdot)$ with $|X^{(r)}(0)|=r$ and $r\to\infty$, 
and with $X^{(r)}(0)/r \to x(0) \in \{-1,1\}$, then in our case it is easy to see that
\beql{eq-fluid-stabil-wp1}
X^{(r)}(rt)/r \stackrel{u.o.c.}{\longrightarrow}  x(t), ~~\mbox{w.p.1},
\eeql
where $x(\cdot)$ is the deterministic Lipschitz function such that $x'(t) = 1-v>0$ when $x(t)<0$, $x'(t) = -v<0$ when $x(t)>0$, 
and therefore $x(t)$ stays at $0$ after it ``hits'' $0$. This implies positive recurrence.

Moreover, since the jump size distribution has finite $\ell$-th moment, $\ell\ge 2$,
using the results of \cite{DM95}, we obtain \eqn{eq-moments-oper}. Since, formally, our model is not within the framework of \cite{DM95},
we provide some details here. In our model the process $X(\cdot)$, with values in $\R$,
corresponds to the queue length process $Q(\cdot)$ in \cite{DM95}, with  values of $Q(t)$ being finite-dimensional vectors with non-negative integer components. In our case the process state norm is $|X(t)|$. 
Property \eqn{eq-fluid-stabil-wp1}, in particular, implies the existence of $T>0$
such that 
\beql{eq-fluid-stabil-weak}
\frac{1}{r} |X^{(r)}(rT)| \stackrel{P}{\longrightarrow} 0.
\eeql
For our model, 
the analog of Proposition 5.1 in \cite{DM95} has the following form. There exists $T>0$ such that,
\beql{eq-moment-to-zero}
\lim_{r\to\infty} \frac{1}{r^\ell} \E |X^{(r)}(rT)|^\ell=0,
\eeql
where $X^{(r)}(\cdot)$ is a sequence of processes defined above in this proof.
Given that we have \eqn{eq-fluid-stabil-weak},
to prove \eqn{eq-moment-to-zero} it suffices to show that the family of
random variables 
$$
\frac{1}{r^\ell} |X^{(r)}(rT)|^\ell
$$
is uniformly integrable. Denote by $U_i$, $i=1,2,\ldots,$ the sequence of i.i.d. random variables, where 
$U_i$ is equal in distribution to a jump size; we have $\E U_i^\ell< \infty$ by \eqn{eq-moment-finite}.
Then the value of $|X^{(r)}(rT)|$ is stochastically dominated by
$$
r + vrT + \sum_{i=1}^{\Pi( rT)} U_i,
$$ 
where $\Pi(t)$ is Poisson process with constant rate $1$.
We see that, to show the uniform integrability of $|X^{(r)}(rT)|^\ell/r^\ell$, it suffices to show that the family 
$$
Z(t) = \frac{[\sum_{s=1}^{\Pi(t)} U_i]^\ell}{t^\ell}
$$
is uniformly integrable for any positive increasing sequence $t\uparrow\infty$; this is true, because $Z(t)\ge 0$ and, as $t\to\infty$, 
we have both $Z(t) \stackrel{P}{\longrightarrow} [\E U_1]^\ell$ and $\E Z(t) \to [\E U_1]^\ell$.
This proves \eqn{eq-moment-to-zero}. 
Since we have \eqn{eq-moment-to-zero}, which is the analog of Proposition 5.1 in \cite{DM95}, 
the rest of the proof of Theorem 5.5 (and thenTheorem 4.1 (i)) in \cite{DM95} applies as is. This, for our model,
proves \eqn{eq-moments-oper}. 
\end{proof}

From now on in the paper, when we consider a distribution $\phi$ and the corresponding distribution $A \phi$, we will denote by $X(\cdot)$ the corresponding single-particle movement process, which defines $A \phi$.

\begin{lem}
\label{lem-A-prop}
Consider any distribution $\phi$, and the corresponding process  $X(\cdot)$, defining $A \phi$. Then
the distribution $A\phi$ is $1/v$-Lipschitz.
\end{lem}

\begin{proof}  
%(i) 
Recall that $\gamma=A \phi$ is nothing else but the stationary distribution of the single-particle movement process $X(\cdot)$, 
corresponding to $\phi$. Consider the stationary regime of this process. Consider any $x$, and any $\Delta>0$. For any fixed initial position $X(0)$ of the particle and any $t>0$, the expected total time in $[0,t]$ that the particle spends
in the interval $(x,x+\Delta]$ is upper bounded by
$$
(t + 1) \Delta/v,
$$
where $(t + 1)$ is the expected number of time intervals {\em between jumps} the particle will have in $[0,t]$, and 
$\Delta/v$ is the maximum time the particle can possibly spend in $(x,x+\Delta]$ during each of those time intervals.
Therefore, considering the stationary regime of the process,
$$
\gamma(x+\Delta) - \gamma(x) \le \frac{[t + 1] \Delta/v}{t}.
$$
This is true for any $t>0$, which implies $\gamma(x+\Delta) - \gamma(x) \le \Delta/v$, i.e. the $1/v$-Lipschitz continuity of $\gamma(\cdot)$.
\end{proof}

\begin{lem}
\label{lem-fp}
A probability distribution function $\phi(\cdot)$ is a traveling wave shape if $\phi = A \phi$.
%For a  probability distribution function $\phi(\cdot)$,  if $\phi = A \phi$ then $\phi$ is a traveling wave shape.
\end{lem}

It can be shown that the converse is also true, and therefore a probability distribution function $\phi(\cdot)$ is a traveling wave shape {\em if 
and only if} $\phi = A \phi$. We do not need the `only if' part for the proof of our main results.

\begin{proof}[Proof of Lemma~\ref{lem-fp}]
 Suppose $\phi = A \phi$. Then, by Lemma~\ref{lem-A-prop}, $\phi$ is Lipschitz. 
Consider any point $x$ where the proper 
derivative $\phi'(x)$ exists. We need to show that \eqn{eq-wave2} holds. 
Consider the single-particle movement process $X(\cdot)$, 
defining $A\phi$.
In the stationary regime, the average rate at which the particle crosses
point $x$ from right to left (``leaves'' $(x,\infty)$, to be precise)  is $v \phi'(x)$; indeed, 
it is easy to see that the expected number of such crossings in interval $[t,t+\Delta t]$ is equal,
up to addition of $o(\Delta t)$ terms, to the probability that at time $t$ the particle is in $(x,x + w \Delta t]$ and it does not jump in $[t,t+\Delta t]$ -- the latter probability itself is equal to $\hat\gamma'(x) w \Delta t + o(\Delta t)$. 
The average rate, in stationary regime,
at which the particle crosses from left to right (``enters'' $(x,\infty)$, to be precise) is the RHS of \eqn{eq-wave2};
indeed, such crossings may occur only when the particle jumps, and the RHS of \eqn{eq-wave2} is the steady-state probability 
that an urge to jump results in a jump over point $x$.
Those right-to-left and left-to-right crossing rates must be equal, thus proving \eqn{eq-wave2}. 
\end{proof}

\section{Traveling wave shapes within a finite frame}
\label{sec-wave-finite-frame}

Let parameters $w \in (0,\infty)$ and $B_L,B_R \in (0,\infty)$ be fixed. 
We will now define operator $A^{(w;B_L,B_R)}$ for a ``finite-frame'' system; this operator is an analog 
of the operator $A$ for the original system. Fixed points of operator $A^{(w;B_L,B_R)}$ will in turn define traveling wave shapes within a finite frame.

\subsection{Operator $A^{(w;B_L,B_R)}$ definition and properties}

Consider a probability distribution function $\gamma=\gamma(\cdot)$ on $\R$, which is actually concentrated on the interval
(``finite frame'') 
$[-B_L,B_R]$, that is $\gamma(-B_L-) =0, \gamma(B_R)=1$. Consider a single particle, whose location $X(t)$ evolves within the 
interval $[-B_L,B_R]$ as follows.
The particle continuously moves left and sometimes makes jumps right. Between jumps,
the particle moves left at the constant speed $w$, unless/until it reaches the left boundary $-B_L$, in which case the particle 
stops at $-B_L$ and stays there until next jump -- this is the ``regulation'' at the left boundary.
The particle get urges to jump as a rate $1$ Poisson process, and actually jumps with probability 
$\bar \eta(X(t),\gamma)$, depending on the current particle location $X(t)$ and distribution $\gamma$.
(Here we {\em cannot} replace $\bar \eta(X(t),\gamma)$ with $\eta(\gamma(X(t)))$, because process $X(\cdot)$ may spend
non-zero time at the left boundary $-B_L$, where $\gamma(\cdot)$ may have a discontinuity.)
The jump sizes are i.i.d. given by the distribution $J(\cdot)$; however, if a jump were to take the particle to the right of $B_R$, the 
particle lands at $B_R$ instead -- this is the ``regulation'' at the right boundary.
Then, the distribution $A^{(w;B_L,B_R)} \gamma$ is defined as the stationary distribution of the process $X(\cdot)$. It is easy to see that process $X(\cdot)$ is positive Harris recurrent: the times when the particle hits $-B_L$ are renewal points and the mean inter-renewal times are clearly finite. Therefore, $X(\cdot)$ indeed has unique stationary distribution. Clearly, $A^{(w;B_L,B_R)} \gamma$ is also concentrated on $[-B_L,B_R]$.

\begin{lem} 
\label{lem-fp-basic-oper} 
If $\hat\gamma=A^{(w;B_L,B_R)} \gamma$, then $\hat\gamma$ satisfies the following conditions:
\beql{eq-lipschitz} 
\mbox{$\hat\gamma(\cdot)$ is $1/w$-Lipschitz in $[-B_L,B_R]$;}
\eeql
\beql{eq-wave-tr-add-oper}
\hat\gamma(-B_L-) =0 < \hat\gamma(-B_L), ~~ \hat\gamma(B_R) = 1;
\end{equation}
almost all points $-B_L < x < B_R$ (w.r.t. Lebesgue measure) are {\em regular}, in that a proper derivative $\hat\gamma'(x)$ exists, and 
at each regular point
\beql{eq-wave-tr-oper}
w \hat\gamma'(x) = \hat\zeta(x) \doteq \int_{-B_L-}^x d\hat\gamma(y) \bar \eta(y;\gamma)\bar J(x-y) ; %     ~~~ -B_L < x < B_R,
\end{equation}
\beql{eq-deriv-right-cont}
\mbox{function $\hat\zeta(x)$ is RCLL  in $[-B_L,  B_R]$, and $\hat\zeta(x-)\ge \hat\zeta(x)$,}
\end{equation}
and
\beql{eq-wave-tr-oper333}
w \frac{d^+}{dx}\hat\gamma(x) = \hat\zeta(x), ~ \forall x\in [-B_L,  B_R),  ~~ \mbox{and}~~
w \frac{d^-}{dx}\hat\gamma(x) = \hat\zeta(x-), ~ \forall x\in (-B_L,  B_R].
\end{equation}
\end{lem}

\begin{proof} 
Recall that $\hat\gamma$ is nothing else but the stationary distribution of the location process $X(\cdot)$ of a single particle,
corresponding to $\gamma$. Consider the stationary regime of this process. Then, the $1/w$-Lipschitz property of $\hat\gamma$
in $(-B_L,B_R]$
is proved in exactly
same way as in the proof of Lemma~\ref{lem-A-prop}, with $w$ replacing $v$. Given that $\hat\gamma$ is right-continuous, it is 
$1/w$-Lipschitz in $[-B_L,B_R]$ as well. Properties \eqn{eq-wave-tr-add-oper} are obvious from the structure of the particle movement process $X(\cdot)$: clearly, in steady-state the process spends zero fraction of time exactly at the right boundary point $B_R$ and non-zero fraction of time at the left boundary point $-B_L$.
 
Since $\hat\gamma(\cdot)$ is Lipschitz in $(-B_L,B_R]$, it is absolutely continuous, and then 
almost all points $-B_L <x<B_R$ are regular, that is the proper 
derivative $\hat\gamma'(x)$ exists.
Consider any regular point $-B_L <x<B_R$. To prove \eqn{eq-wave-tr-oper}, 
essentially the same argument as in the proof of Lemma~\ref{lem-fp} applies.
Indeed, in the stationary regime, the average rate at which the particle crosses
point $x$ from right to left (``leaves'' $(x,B_R]$, to be precise)  is $w \hat\gamma'(x)$ (by the same argument as in the proof of Lemma~\ref{lem-fp}). 
The average rate, in stationary regime,
at which the particle crosses from left to right (``enters'' $(x,B_R]$, to be precise) is $\hat\zeta(x)$;
indeed, such crossings may occur only when the particle jumps, and $\hat\zeta(x)$ is the steady-state probability 
that an urge to jump results in a jump over point $x$.
Those right-to-left and left-to-right crossing rates must be equal, which gives \eqn{eq-wave-tr-oper}. 

Property \eqn{eq-deriv-right-cont} easily follows by directly analyzing the expression for $\hat\zeta(x)$, using the facts that: 
$\bar J(\cdot)$ is non-increasing RCLL; the measure given by $\hat\gamma$ has exactly one atom at $-B_L$ and is absolutely continuous
in $(-B_L,B_R]$. We omit details. Finally, \eqn{eq-wave-tr-oper333} follows from \eqn{eq-deriv-right-cont} and \eqn{eq-wave-tr-oper}.
\end{proof}

Let $\Gamma(w;B_L,B_R)$ denote the set of distribution functions $\gamma$ concentrated on interval $[-B_L,B_R]$, $B_L,B_R>0$, which are
$1/w$-Lipschitz on $[-B_L,B_R]$; here $w>0$, $B_L>0$ and $B_R>0$ are parameters ``attached to'' $\gamma$.
We denote $\Gamma = \cup_{w,B_L,B_R>0} \Gamma(w;B_L,B_R)$.
In particular, $\gamma\in\Gamma$ implies that
 $\gamma(B_R-)=\gamma(B_R)=1$ and $\gamma(B_L-)=0 \le \gamma(B_L)$. So, $\gamma$ may have a single atom, at $-B_L$.
We endow $\Gamma$ with the following natural topology: the convergence $\gamma^{(k)} \stackrel{\Gamma}{\rightarrow} \gamma$ is defined as 
$w^{(k)} \to w$, 
$B_L^{(k)} \to B_L$, $B_R^{(k)} \to B_R$, and $\gamma^{(k)}(y) \to \gamma(y)$ for all $y\in (-B_L,B_R)$. 
In particular, $\gamma^{(k)} \stackrel{\Gamma}{\rightarrow} \gamma$ implies $\gamma^{(k)} \stackrel{J_1}{\rightarrow} \gamma$.
In the special case when $w^{(k)} \equiv w$, 
$B_L^{(k)} \equiv B_L$, $B_R^{(k)} \equiv B_R$, convergence $\gamma^{(k)} \stackrel{\Gamma}{\rightarrow} \gamma$ is equivalent to 
$\gamma^{(k)} \stackrel{u}{\rightarrow} \gamma$.

\begin{lem}%[Map $A^{(w;B_L,B_R)}$ continuity] 
\label{map-monotone}
Map $A^{(w;B_L,B_R)}\gamma$, with $(w;B_L,B_R)$ being the parameters of $\gamma$, 
is continuous in $\gamma \in \Gamma$ (in $\Gamma$-topoogy).
\end{lem}

\begin{proof} 
As $\gamma^{(k)} \stackrel{\Gamma}{\rightarrow} \gamma$, denote the corresponding single-particle processes 
$X_k(\cdot)$ and $X(\cdot)$. Clearly, the sequence $A^{(w^{(k)};B_L^{(k)},B_R^{(k)})}\gamma^{(k)}$ is tight. Consider its 
any subsequential limit $\hat \gamma$, along a subsequence of $k$. Consider each of the processes $X_k(\cdot)$ and $X(\cdot)$ on a 
finite time interval $[0,T]$. These processes can be naturally coupled so that the following property holds w.p.1 along the chosen subsequence. 
If either (a) $X(0)=-B_L$ and $X_k(0) = -B_L^{(k)}$ for all $k$ or (b) $X_k(0) \to X(0) > -B_L$, then 
$$
(X_k(t), ~ 0 \le t \le T) ~~ \stackrel{J_1}{\rightarrow} ~~ (X(t), ~ 0 \le t \le T).
$$
Given this kind of continuity, it is straightforward to see that the limit  $\hat \gamma(\cdot)$ 
of the stationary distributions of $X_k(\cdot)$
must be the stationary distribution $A^{(w;B_L,B_R)}\gamma$ of $X(\cdot)$.
\end{proof}

\subsection{Finite-frame traveling wave shape definition and properties}

Denote by $\mathcal T(w;B_L,B_R)$ the set of those distributions $\gamma$ which are fixed points (FP) 
of $A^{(w;B_L,B_R)}$, i.e. satisfy $\gamma=A^{(w;B_L,B_R)}\gamma$. Any $\gamma\in \mathcal T(w;B_L,B_R)$ we will call
{\em a traveling wave shape for the finite frame $[-B_L,B_R]$ and speed $w$}. If such a traveling wave shape $\gamma$ exists and is unique, i.e. the cardinality of $\mathcal T(w;B_L,B_R)$ is $1$, we will slightly abuse notation by writing $\gamma=\mathcal T(w;B_L,B_R)$.

\begin{lem} 
\label{lem-fp-new}
Let $w,B_L,B_R>0$ be fixed.

(i) $\gamma = \mathcal T(w;B_L,B_R)$ exists and is unique.

(ii) $\gamma = \mathcal T(w;B_L,B_R)$ if and only if it satisfies the following conditions:
\beql{eq-lipschitz-fp} 
\mbox{$\gamma(\cdot)$ is $1/w$-Lipschitz in $[-B_L,B_R]$;}
\eeql
\beql{eq-wave-tr-add-fp}
\gamma(-B_L-) =0 < \gamma(-B_L), ~~ \gamma(B_R) = 1;
\end{equation}
almost all points $-B_L < x < B_R$ (w.r.t. Lebesgue measure) are {\em regular}, in that a proper derivative $\gamma'(x)$ exists, and 
at each regular point
\beql{eq-wave-tr-fp}
w \gamma'(x) = \zeta(x) \doteq \int_{-B_L-}^x d\gamma(y) \bar \eta(y;\gamma)\bar J(x-y) ; 
\end{equation}
\beql{eq-deriv-right-cont-fp}
\mbox{function $\zeta(x)$ is RCLL  in $[-B_L,  B_R]$, and $\zeta(x-)\ge \zeta(x)$,}
\end{equation}
and
\beql{eq-wave-tr-fp333}
w \frac{d^+}{dx}\gamma(x) = \zeta(x), ~ \forall x\in [-B_L,  B_R),  ~~ \mbox{and}~~
w \frac{d^-}{dx}\gamma(x) = \zeta(x-), ~ \forall x\in (-B_L,  B_R];
\end{equation}
\beql{eq-deriv-positive} 
\inf_{[-B_L,  B_R]} \zeta(x) > 0.
\eeql
\end{lem}

\begin{proof}[Proof of Lemma~\ref{lem-fp-new}]
The proof consists of the following three claims.

{\em Claim 1.} $\gamma \in \mathcal T(w;B_L,B_R)$ exists.

{\em Claim 2.} If $\gamma \in \mathcal T(w;B_L,B_R)$ then it satisfies conditions \eqn{eq-lipschitz-fp}-\eqn{eq-deriv-positive}.

{\em Claim 3.} A distribution $\gamma$ concentrated on $[-B_L,B_R]$ and satisfying \eqn{eq-lipschitz-fp}-\eqn{eq-deriv-positive} is unique.

\begin{proof}[Proof of Claim 1] The existence follows from the Brouwer fixed point theorem. (Cf. \cite{Kantorovich-82}, Theorem XVI.5.1, 
for its more general form -- Kakutani theorem.) Indeed, 
by Lemma~\ref{map-monotone}, for a fixed set of positive parameters $w, B_L,B_R$, $A^{(w; B_L,B_R)}\gamma$ continuously (in uniform convergence topology) maps $\Gamma(w;B_L,B_R)$ into itself. Set $\Gamma(w;B_L,B_R)$ is convex and compact. Therefore, a fixed point $\gamma = A^{(w; B_L,B_R)}\gamma$ exists.
%End of Claim 1 proof.
\end{proof}

\begin{proof}[Proof of Claim 2]
Properties \eqn{eq-lipschitz-fp}-\eqn{eq-wave-tr-fp333} follow  from Lemma~\ref{lem-fp-basic-oper}, 
properties
\eqn{eq-lipschitz}-\eqn{eq-wave-tr-oper333}.

Next, observe that $\gamma(x)$ must be strictly increasing in $[-B_L,B_R]$. If not, we could find an interval $[y_1,y_2] \subset [-B_L,B_R)$,
such that $y_1 < y_2$, $\gamma(y_2)=\gamma(y_1)$ (and then $(d^+/dx) \gamma(y_1)=0$), and $\gamma(x) < \gamma(y_1)$ for $x< y_1$; but, by \eqn{eq-wave-tr-fp333} and \eqn{eq-wave-tr-fp}, $(d^+/dx) \gamma(y_1) = \zeta(y_1)/w > 0$, a contradiction.

Let us now prove \eqn{eq-deriv-positive}. First, from \eqn{eq-wave-tr-fp} and the fact that $\gamma(x)$ is strictly increasing
in $[-B_L,B_R]$, observe that $\liminf_{x \uparrow B_R} \zeta(x) > 0$. Therefore, if $\inf_{[-B_L,  B_R]} \zeta(x) = 0$ were to hold, we would 
have a point $y\in [-B_L,B_R)$ such that $\zeta(y)=0$. (Here we use \eqn{eq-deriv-right-cont-fp}.) But $\gamma(x)$ is strictly increasing in a neighborhood of $y$,
which, again by \eqn{eq-wave-tr-fp}, implies $\zeta(y)>0$. This completes the proof of \eqn{eq-deriv-positive} and of Claim 2.
\end{proof}

\begin{proof}[Proof of Claim 3] 
First, note that any $\gamma$ satisfying \eqn{eq-lipschitz-fp}-\eqn{eq-deriv-positive} is such that the following holds.
The inverse $\gamma^{-1}(\nu)$ for $\nu\in (0,1]$ is continuous and Lipschitz, and is strictly increasing in $[\gamma(-B_L),1]$. Moreover, mapping $\nu=\gamma(x)$ gives the one-to-one correspondence between regular points of $\gamma$ in $[-B_L,B_R]$ and regular points of 
$\gamma^{-1}$ in $[\gamma(-B_R),1]$, with the derivatives at the corresponding regular points satisfying
$(d/d\nu) \gamma^{-1}(\nu)= 1/\gamma'(x)$. Furthermore, the latter relation holds for the right derivative for any $x\in [-B_L,B_R)$ and corresponding $\nu=\gamma(x)$, and for the left derivative for any $x\in (-B_L,B_R]$.

The proof of uniqueness is by contradiction. 
Suppose $\gamma_1$ and $\gamma_2$ are two different distributions satisfying conditions \eqn{eq-lipschitz-fp}-\eqn{eq-deriv-positive}.
Let us denote: $\alpha_i(x) \doteq w (d^+/dx) \gamma_i(x)$, $i=1,2$;
then, by the first equation in \eqn{eq-wave-tr-fp333}, $\alpha_i(x)=\zeta_i(x)$.

Suppose, for example, that $\gamma_1(y) > \gamma_2(y)$ for at least one $y\ge 0$. Then let us call
$$
h = \max_\nu [\gamma_2^{-1}(\nu) - \gamma_1^{-1}(\nu)] >0
$$
the ``horizontal distance'' between $\gamma_1$ and $\gamma_2$. It is well-defined, due to the properties of $\gamma_1$ and $\gamma_2$.
Let $\nu_*$ be the minimum $\nu$, at which the horizontal distance is attained. Clearly, $0<\nu_*<1$.
Denote $y_1 = \gamma_1^{-1}(\nu_*) < y_2 = \gamma_2^{-1}(\nu_*)$. Then we must have 
$(d^+/dx) \gamma_1(y_1) \le (d^+/dx) \gamma_2(y_2)$ or, equivalently, $\alpha_1(y_1)\le \alpha_2(y_2)$.
Given the choice of $\nu_*$, $y_1 - \gamma_1^{-1}(\nu) \le y_2 - \gamma_2^{-1}(\nu)$ for all $\nu \in [0,1]$, and then $\zeta_1(y_1) \ge \zeta_2(y_2)$.
Since $\alpha_i(y_i)=\zeta_i(y_i)$, we must have
$$
\alpha_1(y_1)= \alpha_2(y_2) = \zeta_1(y_1) = \zeta_2(y_2). 
$$

Note that if, for example, $\bar J (y) <1$ for all $y>0$, we immediately obtain a contradiction, because in this case $\zeta_1(y_1)> \zeta_2(y_2)$.
The case of general jump size distribution $J$ requires a bit more details, which are as follows.

We must have that $\bar J(y_2-)=1$, i.e. the jump size is at least $y_2$. Indeed, 
by the choice of $\nu_*$, $y_1 - \gamma_1^{-1}(\nu) < y_2 - \gamma_2^{-1}(\nu)$ for all $\nu < \nu_*$. Also, recall that
both inverse functions $\gamma_i^{-1}(\nu)$ are continuous. This means that as we increase $\nu$ in the interval $[0,\nu_*)$,
the interval $[y_1 - \gamma_1^{-1}(\nu),y_2 - \gamma_2^{-1}(\nu)]$ continuously changes (``moves left'') from 
$[y_1 + B_L, y_2 + B_L]$ to (in the limit!) $[0,0]$, while remaining non-zero length. If at least one of points $y \in [0,y_2)$ would be a point of decrease of $\bar J$, we would have $\bar J (y_1 - \gamma_1^{-1}(\nu)) >  \bar J (y_2 - \gamma_2^{-1}(\nu))$ on a non-zero Lebesgue measure subset of $[0,\nu_*)$, which would imply $\zeta_1(y_1) > \zeta_2(y_2)$, which would contradict the definition of $\nu_*$.

Consider now 
$$
\nu^* = \max \{\bar\nu \ge \nu_* ~|~ \gamma_2^{-1}(\nu) - \gamma_1^{-1}(\nu) = h, ~\forall \nu \in [\nu_*,\bar\nu]\}.
$$
We must have $\nu^* < 1$. Denote $y_i^* = \gamma_i^{-1}(\nu^*)$. It is easy to see that $y_2^* = \max \{y \ge y_2 ~|~ \bar J(y-)=1\}$, i.e. $y_2^*$ must be exactly the smallest point of decrease of $\bar J$. Indeed, $\bar J$ cannot have a point of decrease $\hat \nu \in [\nu_*,\nu^*)$. If such point $\hat \nu$ existed, then using the argument analogous to that we used just above to show that $\bar J(y_2-)=1$, we would obtain that $\gamma_1^{-1}(\nu)- \gamma_2^{-1}(\nu)$ would be strictly increasing in a small interval immediately to the right of point $\hat \nu$, which would contradict the definition of the horizontal distance $h$. On the other hand, if $y_2^*$ is {\em not} a point of decrease of $\bar J$, i.e. $\bar J(\hat \nu)=1$ for some $\hat y > y_2^*$,
then again using essentially same argument, we would obtain that $\gamma_1^{-1}(\nu)- \gamma_2^{-1}(\nu)$ would have 
to remain constant in 
a small interval immediately to the right of point $\nu^*$, which would contradict the definition of  $\nu^*$. 

Finally, if $y_2^*$ is the point of decrease of $\bar J$, then, using essentially same argument once again, we obtain that 
$\gamma_1^{-1}(\nu)- \gamma_2^{-1}(\nu)$ must be strictly increasing in a small interval immediately to the right of point $\nu^*$,
which contradicts the definition of the horizontal distance $h$. 
%End of Claim 3 proof.
\end{proof}

The proof of Lemma~\ref{lem-fp-new} is complete.
\end{proof}

\begin{lem}[Continuity, monotonicity, and shift properties] 
\label{lem-monotone}
(i) $\gamma=\mathcal T(w; B_L,B_R)$, as an element of $\Gamma$, is continuous in $(w; B_L,B_R)$. 

(ii) $\gamma=\mathcal T(w; B_L,B_R)$, as a probability distribution, is monotone (in the sense of $\preceq$) in each of the parameters $w$, $B_L$ and
$B_R$, given the other two parameters are fixed. Namely, $\gamma$ is monotone non-increasing in $w$ and $B_L$, and monotone non-decreasing in $B_R$.

(iii) Suppose $B_L$ and $B_R$ are fixed. If
$w\uparrow\infty$ [resp., $w\downarrow 0$], then $\gamma=\mathcal T(w; B_L,B_R)$, as a probability distribution, weakly converges 
to the Dirac distribution concentrated at $-B_L$ [resp., at $B_R$].
Consequently, for any fixed $B_L,B_R$, any fixed $y\in (-B_L,B_R)$, and any fixed $\nu\in (0,1)$, there exists $w>0$ such that 
$\gamma^{-1}(\nu)=y$. In other words, by changing $w$ we can always move a given quantile of $\gamma$ to any point within 
$(-B_L,B_R)$. 

(iv) For any fixed parameters $w$, $B_L$ and $B_R$, and any $c \in \R$, $\gamma=\mathcal T(w; B_L,B_R)$ implies
$\theta_c \gamma=\mathcal T(w; B_L-c,B_R+c)$. In other words, a shift of the finite frame by $c$ results in the same shift
of the corresponding traveling wave shape.
\end{lem}

\begin{proof} 
(i) The continuity follows immediately from the continuity of operator $A^{(w; B_L,B_R)}$ on $\Gamma$ (Lemma~\ref{map-monotone}), and the uniqueness of $\gamma=\mathcal T(w; B_L,B_R)$ (Lemma~\ref{lem-fp-new}(i)).

(ii) The monotonicity is proved by contradiction, using the horizontal distance, almost exactly same way as in the proof of uniqueness (Claim 3) in Lemma~\ref{lem-fp-new}.

(iii) The case $w\uparrow\infty$ is obvious, because $\gamma$ is $1/w$-Lipschitz in $[B_L,B_R]$ and $\gamma(B_R)=1$.
Suppose $w\downarrow 0$. Then, $\gamma=\mathcal T(w; B_L,B_R)$ is monotone non-increasing in $w$. It suffices to show that, for each 
$y< B_R$, $\lim \gamma(y) =0$. Indeed, fix $y<B_R$ and any $y_1\in (y,B_R)$. Clearly, $\lim \gamma(y_1) < \nu_1 < 1$. Fix small $\epsilon>0$, so that $y < y_1-\epsilon < y_1$. Then, uniformly in all sufficiently small $w>0$, the corresponding particle process $X(\cdot)$ is such that 
when it is to the left of point $y_1$, it will jump with the rate at least $\eta(\nu_1)$, while it moves left at the small speed $w$; this in turn 
easily implies that the steady-state probability of the particle being to the left of $y_1-\epsilon$ vanishes.

(iv) Obvious from the definition of a finite-frame traveling wave shape.
\end{proof}

\section{Proofs of the main results, Theorems~\ref{thm-existence} and \ref{thm-main}}
\label{sec-main-proofs}

\subsection{Proof of Theorem~\ref{thm-existence}(i): Existence}
\label{sec-main-proofs-existence}

The high level outline of the proof is as follows. We obtain a traveling wave shape $\phi$ as a limit of finite-frame traveling wave shapes $\gamma_B$, as $B\to\infty$, for the finite frame $[-B,B]$ and with the speed $w_B$ chosen so that the median of $\gamma_B$ is at $0$. The first step, Lemma~\ref{lem111}, shows that, necessarily, $w_B \to v$. The next step is to show, in Lemma~\ref{lem-gamma_B-tight}, that the family of distributions $\{\gamma_B\}$ is tight. The proof of Lemma~\ref{lem-gamma_B-tight} involves considering the single-particle movement processes 
$X_B(\cdot)$, corresponding to  $\gamma_B$, and their space/time-rescaled versions $X_{(C)}(\cdot)=X_B(C \cdot)/C$, with $C=C(B) \to \infty$
as $B\to\infty$. (Supplementary Lemma~\ref{lem-conv-determ}, used in the proof of Lemma~\ref{lem-gamma_B-tight}, 
is given in Section~\ref{sec-suppl-fact}.) 
After the tightness of $\{\gamma_B\}$ is established, Lemma~\ref{lem5} completes the
proof of Theorem~\ref{thm-existence}(i) by demonstrating that any subsequential limit $\phi$ of $\gamma_B$ must be a traveling wave shape.
At the beginning of the proofs of Lemmas~\ref{lem111},  \ref{lem-gamma_B-tight} and  \ref{lem-conv-determ}
we give intuition and outlines specifically for those proofs.

We now proceed with the formal proof.
Let us consider the finite-frame  system with $B_L= B_R=B$, i.e. within frame $[-B,B]$. Let $B\uparrow\infty$ along some fixed sequence. For each $B$, we will choose the speed $w_B$ such that the corresponding $\gamma$, which we denote by $\gamma_{B}$,
has its median exactly at $0$, i.e. $\gamma_{B}(0)=1/2$ or, equivalently, $\gamma_B^{-1}(1/2)=0$. (It does not have to be the median, 
we could fix any $\nu_0\in (0,1)$ and choose the speed $w_B$ such that $\gamma_B^{-1}(\nu_0)=0$.) 
Such $w_B$ and $\gamma_{B}$ exist by Lemma~\ref{lem-monotone}(iii), and $0 < w_B < \infty$. 
The corresponding single-particle movement process $X(\cdot)$ we will denote by $X_B(\cdot)$.

\begin{lem}
\label{lem111}
Necessarily, $\lim_{B\to\infty} w_B =v$.
\end{lem}

\begin{proof} 
The proof is by contradiction. The basic intuition is as follows. Suppose, 
for example, that $\limsup_{B\to\infty} w_B > v$. Then there exists a constant $w_\infty > v$ such that
$w_B > w_\infty > v$ along some subsequence of $B$. Then we show that the distributions $\gamma_B$ would have to stay tight around the left end $-B$ of the frame. This is, informally speaking, due to the fact that, if in steady-state $X_B(\cdot)$ would stay away 
from $-B$ most of the time, then its average drift would be upper bounded by $v - w_\infty < 0$, while it must be $0$ for any $B$. 
But, if $\gamma_B$ stays tight around $-B$, then as $B\to\infty$, the median of $\gamma_B$ cannot stay at $0$, which contradicts the 
$\gamma_B$ definition. The contradiction in the case $\liminf_{B\to\infty} w_B < v$ is obtained similarly -- the would have
distributions $\gamma_B$  staying tight around the right end $B$ of the frame.

We proceed with the formal proof. Suppose $w_B$ does not converge to $v$.
Consider, first, the case\\ $\limsup_{B\to\infty} w_B > v$. Then, we can choose a constant $w_\infty > v$ and
a subsequence of $B$, along which $w_B > w_\infty > v$. For all large $B$, consider $\gamma_B^{up}=\mathcal T(w_\infty;B,B)$ corresponding to speed $w_\infty$;
also consider the corresponding single-particle process $X_B^{up}(\cdot)$.
By monotonicity in $w$ (Lemma~\ref{lem-monotone}(ii)), $\gamma_B \preceq \gamma_B^{up}$ for all large $B$. By monotonicity in the right end of the frame (Lemma~\ref{lem-monotone}(ii)), and the shift property (Lemma~\ref{lem-monotone}(iv)), 
we observe that the distribution $\theta_B \gamma_B^{up}$ (i.e., the distribution $\gamma_B^{up}$ shifted right by $B$, and therefore concentrated on $[0,2B]$) will be monotone non-decreasing (in the sense of $\preceq$) in $B$. Then $\theta_B \gamma_B^{up}$ must uniformly converge, $\theta_B \gamma_B^{up} \stackrel{u}{\rightarrow} \gamma_{\infty}$
 to some distribution $\gamma_{\infty}$, concentrated on $[0,\infty)$, and the function $\gamma_{\infty}$ is  $1/w_\infty$-Lipschitz
 on  $[0,\infty)$. (This follows from the fact that all functions $\theta_B \gamma_B^{up}$ are uniformly $1/w_\infty$-Lipschitz in $[0,\infty)$.) Moreover, 
\beql{eq-atom-at-0}
\gamma_{\infty}(0)>0.
\eeql
Indeed the particle process $X_B^{up}(\cdot)$ is such that the steady-state average drift of the 
particle is upper bounded by $-(1-\gamma_B^{up}(-B))w_\infty + v$; if \eqn{eq-atom-at-0} would not hold, then, as $B\to\infty$, 
this upper bound would converge to $-(1-\gamma_{\infty}(0))w_\infty + v =-w_\infty + v <0$, which is impossible.

To obtain a contradiction, it remains to show that the distribution $\gamma_{\infty}$ is proper, i.e. $\lim_{y\to\infty} \gamma_{\infty}(y) = 1$.
Indeed, this would imply that the median of 
$\gamma_B^{up}$, and then of $\gamma_B$ as well, will eventually become less than $0$, which contradicts the 
definition of $\gamma_B$.

Let us show that the distribution $\gamma_{\infty}$ is proper. Suppose not, that is $\lim_{y\to\infty} \gamma_{\infty}(y) = \nu^*_\infty < 1$.
Denote $\nu^* = \eta^{-1}(w_\infty)$. 

The case $\nu^*_\infty > \nu^*$ is impossible. Indeed, in this case, pick a small $\epsilon>0$ and any fixed $H>0$ such that $\gamma_{\infty}(H)  > \nu^* + \epsilon$. Then, uniformly in all large $B$, when the 
particle process $X_B^{up}(\cdot)+B$ (``living'' in $[0,2B]$)
 is to the right of $H$, then the particle jumps right at most at the rate $\eta(\nu^* + \epsilon) < w_\infty$. We also have the following: when the process $X_B^{up}(\cdot)+B$ ``jumps over point $H$,'' then
uniformly in all $B$ and on the jump starting point in $[0,H]$, the distribution of the ``overshoot-over-$H$'' distance $V$
is uniformly stochastically upper bounded by some proper distribution with finite mean. (This is obvious if jump size distribution $J(\cdot)$
has finite support. Otherwise, $\pr\{V > y\} \le \bar J(y) / \bar J(H)$.)
Using these properties, we can construct a uniform in $B$ proper stochastic upper bound on the stationary distribution of $X_B^{up}(\cdot)+B$,
which is $\theta_B \gamma_B^{up}$. This would contradict the fact that the distribution $\gamma_{\infty}$ is not proper.

Consider the case $\nu^*_\infty \le \nu^*$. Consider the particle process $X_\infty(\cdot)$, taking values in $\R_+=[0,\infty)$,
and defined exactly as the process defining 
the (finite-frame) operator value $A^{(w;B_L,B_R)} \gamma$, except: $X_\infty(\cdot)$ corresponds to speed $w_\infty$ and distribution $\gamma_\infty$ on $\R_+=[0,\infty)$ (which is not necessarily proper); the lower boundary is $0$ (at which regulation does occur); there is no upper boundary (and no regulation from above -- no forward jump is ever ``truncated'').
Process $X_\infty(\cdot)$ is regenerative, with regenerations occurring when the particle hits $0$.
Given that for any finite $y\ge 0$, $\eta(\gamma_\infty(y)) > w_\infty$, the process cannot possibly be positive recurrent.
(Indeed, even if the rate of jumps is constant at $w_\infty$, i.e. lower than it actually is, 
process $X_\infty(\cdot)$ models workload of a queue, where the workload is depleted at constant rate $w_\infty$ and new workload arrives 
as a random process with average rate exactly $w_\infty$. Therefore, $X_\infty(\cdot)$ is stochastically lower bounded by a  
non-positive recurrent process.)
Consider two sub-cases: (a) $X_\infty(\cdot)$ is transient (which is certainly the case if $\nu^*_\infty < \nu^*$,
but logically also possible if $\nu^*_\infty = \nu^*$;) (b) $X_\infty(\cdot)$ is null-recurrent (only possible when $\nu^*_\infty = \nu^*$).

Observe that the processes $X_B^{up}(\cdot)+B$ for all $B$ and the process $X_\infty(\cdot)$, all starting at $0$, can be naturally coupled so that,
as $B\to\infty$, the trajectory of $X_B^{up}(\cdot)+B$ converges (in Skorohod topology)
to that of $X_\infty(\cdot)$ w.p.1. 

Sub-case (a). If $X_\infty(\cdot)$ is transient, then for some fixed $\delta>0$, and any $H>0$, for all sufficiently large $B$, the probability
that $X_B^{up}(\cdot)+B$ will up-cross level $H$ in a regeneration cycle is at least $\delta$. This would imply that the mean duration of the regeneration cycle of $X_B^{up}(\cdot)+B$ has to go to $\infty$ as $B\to\infty$. (Because, the time to return to $0$ from a point $y\ge H$ is at least 
$H/(w_\infty/2)$.) This would imply that the steady-state probability $\theta_B \gamma_B^{up}(0)$
of $X_B^{up}(\cdot)+B$ being at $0$ will vanish. Indeed, the mean time the process spends at $0$ within one regeneration cycle is uniformly upper bounded by an exponential random variable 
with mean $1/v$ (because, the rate at which the particle jumps forward when it sits at the left boundary point is at least $\int_0^1 \eta(\nu) d\nu = v$). But, vanishing $\theta_B \gamma_B^{up}(0)$ would mean $\gamma_\infty(0)=0$, which contradicts \eqn{eq-atom-at-0}.

Sub-case (b). If $X_\infty(\cdot)$ is null-recurrent, then the regeneration cycle of $X_\infty(\cdot)$ is finite w.p.1, and has infinite mean. 
Recall the natural coupling of $X_B^{up}(\cdot)+B$ and $X_\infty(\cdot)$. Given this coupling, for any fixed $H>0$, on the event that $X_\infty(\cdot)$
does not up-cross $H$ in a regeneration cycle, the regeneration cycle length of $X_B^{up}(\cdot)+B$ converges to that of $X_\infty(\cdot)$ almost surely. Using Fatou's lemma, we see that, as $B\to\infty$, $\liminf$ of the mean regeneration cycle duration of $X_B^{up}(\cdot)+B$  is at least that of $X_\infty(\cdot)$, which is infinity. This, again, leads to vanishing $\theta_B \gamma_B^{up}(0)$ and a contradiction with \eqn{eq-atom-at-0}.

This completes the proof that $\limsup_{B\to\infty} w_B > v$ is impossible. 

The contradiction to $\liminf_{B\to\infty} w_B < v$ is obtained similarly. In this case, we can construct lower bounds on $\gamma_B$. 
Condition \eqn{eq-atom-at-0} is replaced by the condition that the steady-state rate of the particle hitting the right boundary $B$ converges down to a positive number. Regeneration time points are when the particle hits the right boundary $B$. We omit details.
\end{proof}

\begin{lem}
\label{lem-gamma_B-tight} 
The family of distributions $\{\gamma_B\}$ is tight. 
\end{lem}

\begin{proof} 
The proof is by contradiction. Its basic outline is as follows.
If $\{\gamma_B\}$ is not tight, we can find a subsequence of $B$, and a corresponding sequence of scaling factors
$C=C(B)\uparrow\infty$ such that rescaled distributions $\psi_{(C)}(y)= \gamma_B(Cy)$ are such that, say, some $\nu$-quantile ($\nu > 1/2)$ of $\psi_{(C)}$ remains at point $1$. Each $\psi_{(C)}$ is a finite-frame traveling shape; its frame is $[-B/C,B/C]$, the 
corresponding single-particle process is $X_{(C)}(t)=X_B(Ct)/C$, and the speed is still $w_B$.
We show that the family of distributions $\{\psi_{(C)}\}$ is tight. This is done by considering a quadratic Lyapunov function 
for $X_{(C)}(t)$ and its steady-state drift, and showing that $\E |X_{(C)}(\infty)|$ remains uniformly bounded. Given the tightness of the sequence $\{\psi_{(C)}\}$ we consider its subsequence, converging to some proper distribution $\psi_{(\infty)}$. 
Using Lemma~\ref{lem-conv-determ} we show that $\psi_{(\infty)}$ must be concentrated on at most two atoms $z_1\le 0$
and $z_2 \ge 0$. The case of one atom, $z_1=z_2$, is impossible because it leads to contradiction with the $\nu$-quantile ($\nu > 1/2)$ of $\psi_{(C)}$ remaining at point $1$. Then we must have either $z_1<0$ or $z_2>0$. In this case we again use a quadratic Lyapunov function for $X_{(C)}(t)$ to show that, for large $C$, its steady-state drift would be negative; roughly because the probability of 
$X_{(C)}(\infty)\ge z_2$ is about $1-\nu$ and the probability of $X_{(C)}(\infty)\le z_1$, and the drift of the quadratic Lyapunov function
is strictly negative when $X_{(C)}(\infty)\ge z_2 > 0$ and non-positive when $X_{(C)}(\infty)\le z_1\le 0$. 

We proceed with the formal proof by contradiction.
Suppose $\{\gamma_B\}$ is not tight. Then, there exists fixed $\nu>1/2$ and a subsequence $B\uparrow \infty$, along which $C=C(B)=\max\{\gamma_B^{-1}(\nu), - \gamma_B^{-1}(1-\nu)\} \uparrow \infty$. 
WLOG we can choose $\nu$ sufficiently close to $1$, so that
\beql{eq-nu-choice}
\eta(1-\nu) > v > \eta(\nu).
\eeql

Note that, for each $B$,
either $\gamma_B^{-1}(\nu)=C$ or $\gamma_B^{-1}(1-\nu)=-C$; and also,
$\gamma_B(C)-\gamma_B(-C) \ge \nu-(1-\nu)=2\nu-1$.

For each $B$ (with the corresponding $C$), consider the particle process with space compressed by $C$ and time sped up by $C$, i.e. the process
$X_{(C)}(t)= X_B(Ct)/C$. For each $B$, the process lives within the frame $[-C_0,C_0]$, where $C_0=B/C$. It is possible that $C_0\to\infty$
as $B\to\infty$. Denote by $\psi_{(C)}(y)=\gamma_B(Cy)$ the corresponding scaled distribution. 
From now on in this proof, when we say `for a given $B$' or `for a given $C$' we mean `for a given pair $(B,C)$ of $B$ and $C$, corresponding to each other.

Note that, for each $C$, 
either $\psi_{(C)}^{-1}(\nu)=1$ or $\psi_{(C)}^{-1}(1-\nu)=-1$; also, $\psi_{(C)}(1)-\psi_{(C)}(-1) \ge 2\nu-1$.

Let us prove that
\beql{eq-gamma_C-tight}
\mbox{the family of distributions $\{\gamma_C\}$ is tight.}
\eeql
This is trivially true if $C_0$ remains bounded as $B\to\infty$. Therefore, it suffices to prove \eqn{eq-gamma_C-tight} 
for a subsequence of $B$ such that $C_0=B/C \to \infty$. Let us consider such a subsequence.

Let $g(y)=y^2/2$. Denote by $P_{(C)}^t(y,H)$ the transition function of the process $X_{(C)}(t)$.
 According to our definitions, for any fixed $C$, $\psi_{(C)} P_{(C)}^t g = \psi_{(C)} g$ for all $t\ge 0$; here we view the transition function $P_{(C)}^t$ and the distribution $\psi_{(C)}$ as operators (see Section ~\ref{sec-notation}).

Note that $\psi_{(C)} = \mathcal T(w_B; C_0,C_0)$; and that $X_{(C)}(\cdot)$ is the corresponding (finite-frame) particle process.
Denote by $\hat X_{(C)}(\cdot)$ a single particle process, corresponding (like $X_{(C)}(\cdot)$) 
to the distribution function $\psi_{(C)}$ and speed $w_B$, but different from $X_{(C)}(\cdot)$ in that it evolves 
in $(-\infty, \infty)$, i.e. no regulation keeping the process within finite frame $[-C_0,C_0]$ is applied.
The particle continues to move left at speed $w_B$ even when it is at or to the left of $-C_0$, in which case its jump rate
remains the same as if it would be at point $-C_0$, namely $\bar \eta(-C_0, \psi_{(C)})$. 
The particle also can ``jump over'' point $C_0$; when it is at or to the right of $C_0$, it moves left at speed $w_B$ and its jump rate
is as it would be at point $C_0$, namely $\bar \eta(C_0, \psi_{(C)}) = \eta(\psi_{(C)}(C_0)) = \eta(1) =0$.
Denote by $\hat P_{(C)}^t(y,H)$ the transition function of the process $\hat X_{(C)}(t)$. We claim that, for any sufficiently large fixed $C$, 
uniformly in $y\in [-C_0,C_0]$,
\beql{eq-trans-bound}
[P_{(C)}^t g] (y) - [\hat P_{(C)}^t g] (y) \le o(t),
\eeql
where $o(t)$ is a positive function (which may depend on $C$), such that $o(t)/t\to 0$ as $t\to 0$. The proof of \eqn{eq-trans-bound} is given in Section~\ref{sec-proof-append}.

It is easy to check directly that,
uniformly in $y\in [-C_0,C_0]$, 
$$ %\beql{eq-upper-generator-def}
\lim_{t\downarrow 0} \frac{1}{t} [\hat P_{(C)}^t - I] g (y) = \hat G_{(C)} g (y) \doteq - w_B g'(y) + C \eta(\psi_{(C)}(y))\int_0^\infty dJ(C\xi) [g(y+\xi)-g(y)],
$$ %\eeql
and we have
$$ %\beql{eq-upper-generator-bound1}
\hat G_{(C)} g (y) = - w_B y + C \eta(\psi_{(C)}(y))\int_0^\infty dJ(\zeta) [y \zeta/C +(\zeta/C)^2/2]
= 
$$
$$
- w_B y + \eta(\psi_{(C)}(y))[y + (1/2)m^{(2)}/C],
$$ %\eeql
and then
\beql{eq-upper-generator-bound}
\hat G_{(C)} g (y) \le y (-w_B+\eta(\psi_{(C)}(y))) + m^{(2)}/(2C).
\eeql
Note that, by our construction and \eqn{eq-nu-choice},
for all sufficiently large $C$ (recall that $w_B \to v$), we have, for some fixed $\epsilon>0$,
\beql{eq-drift2}
-w_B+\eta(\psi_C(y)) \ge \epsilon, ~ y\le -1, ~~~\mbox{and}~~~ -w_B+\eta(\psi_C(y)) \le -\epsilon, ~ y\ge 1.
\eeql

From \eqn{eq-trans-bound} we can see that, informally speaking, operator $\hat G_{(C)}$ is an upper bound on the generator $G_{(C)}$ of process $X_{(C)}(\cdot)$, when these operators are applied to function $g$. This observation is only informal, because function $g$ is {\em not} even within the domain of  $G_{(C)}$. Formally, using \eqn{eq-trans-bound} and the fact that distribution $\psi_{(C)}$ is
concentrated on $[-C_0,C_0]$, 
we can write: 
\beql{eq-key-generator-ineq}
0 = \lim_{t\downarrow 0} \frac{1}{t} \psi_{(C)} [P_{(C)}^t - I] g \le \lim_{t\downarrow 0} \frac{1}{t} \psi_{(C)} [\hat P_{(C)}^t - I] g
= \psi_{(C)} \hat G_{(C)} g, 
\eeql
and observe that
$$
\psi_{(C)} \hat G_{(C)} g \le -\epsilon \psi_{(C)} h_1 + C_2,
$$
where $h_1(y)\doteq |y| I(|y|\ge 1)$, and
$C_2$ does not depend on (sufficiently large) $C$.
We obtain that, uniformly in large $C$,
$$
\psi_{(C)} h_1 \le C_2/\epsilon.
$$
This implies the tightness of $\{\psi_{(C)}\}$,  thus completing the proof of claim \eqn{eq-gamma_C-tight}.

Given the tightness of $\{\psi_{(C)}\}$, consider a subsequence of $C$ along which 
$\psi_{(C)} \stackrel{w}{\rightarrow} \psi_{(\infty)}$, where $\psi_{(\infty)}$ is a proper distribution, not necessarily continuous. Recall that $w_B \to v$. Then, we can apply Lemma~\ref{lem-conv-determ}.
We obtain that the entire distribution $\psi_{(\infty)}$ must be concentrated at either (a) single point $z\in[-1,1]$ or (b)
on a non-zero length segment $[z_1,z_2]$, where $z_1 = \min\{y~|~\psi_{(\infty)}(y) = \nu_v\}$, $z_2 = \sup\{y~|~\psi_{(\infty)}(y) = \nu_v\}$,
$\nu_v \doteq \eta^{-1}(v)$.  

Case (a) is impossible, because the property that either $\psi_{(C)}^{-1}(\nu)=1$ or $\psi_{(C)}^{-1}(1-\nu)=-1$ cannot hold for all large $C$. 

Let us show that case (b) is also impossible. Suppose, it holds. Then, $z_1 \le 0 \le z_2$, because otherwise the median of $\psi_{(C)}$ could not stay at $0$ for all $C$. We also know that (since $\psi_{(\infty)}(y)=\nu_v \in [z_1,z_2)$) the distribution $\psi_{(\infty)}$ consists of two atoms at $z_1$ and $z_2$, with positive weights $\nu_v$ and $1-\nu_v$, respectively. So, either $z_1<0$ or 
$z_2>0$. Consider, for concreteness, the case $z_2>0$. (The treatment of the case $z_1<0$ is analogous.) Fix small $\epsilon>0$ and, for each $C$, consider the points
$$
z_{1,\epsilon} = z_{1,\epsilon}(C) \doteq \psi_{(C)}^{-1}(\nu_v -\epsilon), ~~ z_{2,\epsilon} = z_{2,\epsilon}(C) \doteq \psi_{(C)}^{-1}(\nu_v +\epsilon).
$$
Given the form of the limiting distribution $\psi_{(\infty)}$, as $C\to\infty$, we have 
$$
z_{1,\epsilon} \to z_1, ~ z_{2,\epsilon} \to z_2.
$$
Let $\epsilon_1= \eta(\nu_v-\epsilon)- \eta(\nu_v+\epsilon)$. Obviously, $\epsilon_1 \downarrow 0$ as $\epsilon\downarrow 0$.

For any fixed $\epsilon$ (and corresponding $\epsilon_1$), and arbitrarily small $\delta>0$, the following holds  
for all sufficiently  large $C$:
\beql{eq-upper-generator-bound2}
\hat G_{(C)} g (y) \le y (-w_B+\eta(\psi_{(C)}(y))) + \delta, ~~\mbox{(from \eqn{eq-upper-generator-bound})}.
\eeql
We will now estimate $\psi_{(C)} \hat G_{(C)} g$ by considering intervals $(-\infty, z_{1,\epsilon}]$, $(z_{1,\epsilon}, z_{2,\epsilon}]$,  
$(z_{2,\epsilon}, \infty)$, and using: the $\psi_{(C)}$-measures of those intervals; estimates of values of $-w_B+\eta(\psi_{(C)}(y))$ over those intervals; estimates of $|z_{1,\epsilon} - z_1|$ and $|z_{2,\epsilon} - z_2|$. The values and estimates listed just above are:
$$
|z_{1,\epsilon} - z_1| < \delta, ~~ |z_{2,\epsilon} - z_2| < \delta,
$$
$$
-w_B+\eta(\psi_{(C)}(y)) > \epsilon_1, ~ y \le z_{1,\epsilon},
$$
$$
-w_B+\eta(\psi_{(C)}(y)) < \epsilon_1, ~ y \ge z_{2,\epsilon},
$$
$$
|-w_B+\eta(\psi_{(C)}(y))| \le \epsilon_1, ~ z_{1,\epsilon}\le y \le z_{2,\epsilon},
$$
$$
\psi_{(C)}(z_{1,\epsilon}) = \nu_v-\epsilon , ~ \psi_{(C)}(z_{2,\epsilon}) - \psi_{(C)}(z_{1,\epsilon}) = 2\epsilon, ~ 1- \psi_{(C)}(z_{2,\epsilon}) = 1-\nu_v - \epsilon.
$$
Then, from \eqn{eq-upper-generator-bound2}, 
$$
\psi_{(C)} \hat G_{(C)} g \le (\nu_v-\epsilon) (z_1+\delta) \epsilon_1 + 2\epsilon [\max\{|z_1|,z_2\} + 2\delta] \epsilon_1 + (1-\nu_v + \epsilon) (z_2-\delta) (-\epsilon_1) +\delta
$$
$$
\le (\nu_v-\epsilon) \delta \epsilon_1 + 2\epsilon [\max\{|z_1|,z_2\} + 2\delta] \epsilon_1 + (1-\nu_v + \epsilon) (z_2-\delta) (-\epsilon_1) +\delta,
$$
where the last inequality is because $z_1 \le 0$.
If we pick sufficiently small $\epsilon>0$ (with the corresponding, also small, $\epsilon_1>0$), and then choose sufficiently small $\delta$, then the RHS of in the last display is negative. This implies $\psi_{(C)} \hat G_{(C)} g < 0$,
which is a contradiction with \eqn{eq-key-generator-ineq}. Therefore, case (b) is impossible.

Thus, the assumption that the family of distributions $\{\gamma_B\}$ is not tight leads to a contradiction. This completes the proof of the lemma.
\end{proof}

The following lemma completes the
proof of Theorem~\ref{thm-existence}(i).

\begin{lem}
\label{lem5} 
Consider a subsequence of $B$. (We know that  $w_B \to v$.)
Then there exists a further subsequence along which $\gamma_{B}$ converges to a proper distribution $\phi$ on $(-\infty,\infty)$,
which is a traveling wave shape with the median at $0$, $\phi(0)=1/2$.
\end{lem}

\begin{proof} 
Using the tightness of $\{\gamma_B\}$, we can find a further subsequence of $B$, along which $\gamma_B \stackrel{w}{\rightarrow} \phi$;
but, since the functions $\gamma_B$ are uniformly Lipschitz in $(-B,\infty)$, we conclude that the function $\phi$ is $1/v$-Lipschitz,
and then $\gamma_B \stackrel{u}{\rightarrow} \phi$.
Also, clearly $\phi(0)=1/2$.

If $P_{B}^t(y,H)$ is the transition function of the process $X_B(\cdot)$, then it is easy to see that, for each $y$ and $t\ge 0$, 
$P_{B}^t(y,\cdot) \stackrel{w}{\rightarrow} P_\infty^t(y,\cdot)$, where $P_\infty^t(y,\cdot)$ is the transition function of the particle process
$X_\infty(\cdot)$, corresponding to distribution $\phi$ and speed $v$. Moreover, 
$P_\infty^t(y,\cdot)$  is continuous in $y$, i.e. process $X_\infty(\cdot)$ is Feller continuous. Then (see, e.g., \cite{Liptser_Shiryaev}), the limit $\phi$ of the stationary distributions $\gamma_B$ of $X_B(\cdot)$ is the stationary distribution of the limiting process $X_\infty(\cdot)$. This means that $\phi$ is in fact the fixed point $\phi=A \phi$, and therefore, by Lemma~\ref{lem-fp}, $\phi$ is a traveling wave shape.
\end{proof}

\subsection{Proof of \eqn{eq-trans-bound}}
\label{sec-proof-append}

To prove \eqn{eq-trans-bound} we first show that, for a small $t$, the contribution into the expectations 
$[P_{(C)}^t g] (y)$ and $[\hat P_{(C)}^t g] (y)$ of the events involving two or more 
particle jump urges in $[0,t]$ is $o(t)$. 
Indeed, consider for example $[\hat P_{(C)}^t g] (y)$ which is equal to \\
$
\E [(\hat X_{(C)}(t))^2/2 | \hat X_{(C)}(0)=y].
$
We need to show that
\beql{eq-two-jumps}
\E [(\hat X_{(C)}(t))^2 | \hat X_{(C)}(0)=y; \mbox{two or more jump urges in $[0,t]$}] \le o(t),
\eeql
for a fixed $C$. 
Let $U_i, ~i=1,2,\ldots$, be an i.i.d. sequence of random variables with the distribution $J(\cdot)$ (of a jump size); let 
$\Pi_t$ be an independent from sequence Poisson random variable with mean $t$; denote $\tilde \Pi_t = \Pi_t I\{\Pi_t \ge 2\}$.
Note that $\E \tilde \Pi_t  = t - te^{-t} = o(t)$, $\E \tilde \Pi_t^2  = t + t^2 - te^{-t} = o(t)$. Denote
$$
Z = [\sum_{i=1}^{\Pi_t} U_i] I\{\Pi_t \ge 2\} = \sum_{i=1}^{\tilde \Pi_t} U_i,
$$
and note that
$$
\E Z \le \E U_1 \E \tilde \Pi_t = o(t), ~~ \E Z^2 = [\E \tilde \Pi_t^2 - \E \tilde \Pi_t] (\E U_1)^2 + \E \tilde \Pi_t \E U_1^2 = o(t).
$$
Finally, for bounded values of $t$, $t \in [0,a]$, 
$$
\E [(\hat X_{(C)}(t))^2 | \hat X_{(C)}(0)=y; \mbox{two or more jump urges in $[0,t]$}] \le 
$$
$$
\E [C_0+w_B a + \sum_{i=1}^{\Pi_t} U_i]^2 I\{\Pi_t \ge 2\} \le
$$
$$
C_{11} \E I\{\Pi_t \ge 2\} + C_{12} \E Z + C_{13} \E Z^2 = o(t),
$$
where constants $C_{11}, C_{12}, C_{13}$ do not depend on $y\in [-C_0,C_0]$. This proves \eqn{eq-two-jumps}.

Continue with the proof of \eqn{eq-trans-bound}. Consider the particle located at $y\in [-C_0,C_0]$ at time $0$.
In the interval $[0,t]$ we couple $X_{(C)}(\cdot)$ and $\hat X_{(C)}(\cdot)$ in the natural way so that,
if the particle has at most one jump urge in this interval, then the particle jump times 
and corresponding jump sizes coincide. Then, if the particle (starting at $y$) has no jumps in $[0,t]$,  $X^2_{(C)}(t) \le \hat X^2_{(C)}(t)$ holds. Consider now the case when the particle has exactly one jump in $[0,t]$, at some time $s$.
 Note that $X_{(C)}(s) \ne \hat X_{(C)}(s)$, specifically
$X_{(C)}(s) > \hat X_{(C)}(s)$, is only possible if the particle ``hits'' point $-C_0$ in $[0,s]$. Also note that, if the latter happens, 
$0 \le \Delta \doteq X_{(C)}(s) - \hat X_{(C)}(s) \le t w_B$ and then, using the relation $g(y+\delta) - g(y) = y\delta + \delta^2/2$, we can write 
$$
\E [g(X_{(C)}(s+)) - g(\hat X_{(C)}(s+))] I\{X_{(C)}(s+) > \hat X_{(C)}(s+)\}
\le 
$$
\beql{eq-bound-diff}
\Delta\int_{0}^{2 C_0} dJ(x) (x-C_0) + \Delta \bar J (2C_0) + o(t).
\eeql
The integral in %the RHS of 
\eqn{eq-bound-diff} is the upper bound on the expectation of the 
linear term of the difference $g(X_{(C)}(s+)) - g(\hat X_{(C)}(s+))$ (the ``$y\delta$'') for the jump sizes $x \in [0,2C_0]$; for such jump sizes,
necessarily, $X_{(C)}(s+) > \hat X_{(C)}(s+)$. The $\Delta \bar J (2C_0)$ is the upper bound on the expected linear term of the difference $g(X_{(C)}(s+)) - g(\hat X_{(C)}(s+))$ for jump sizes $x > 2C_0$; in this case $X_{(C)}(s+) = C_0 \ge \hat X_{(C)}(s+)$, 
and we can use the same bound as for jump size $x \le C_0$. We can write
$$
\mbox{%RHS of 
\eqn{eq-bound-diff}} \le 
- \Delta\int_0^{C_0} dJ(x) (C_0-x) + \Delta\int_{C_0}^{2 C_0} dJ(x) (x-C_0) + \Delta \bar J (2C_0) + o(t) = 
$$
$$
[-C_5+C_6]\Delta + o(t),
$$ 
where 
$$
C_5 = C_0 J(C_0) - \int_0^{C_0} dJ(x) x,
$$
$$
C_6 = \int_{C_0}^{2 C_0} dJ(x) x - C_0 [J(2C_0)-J(C_0)] + \bar J (2C_0).
$$
Note that for all sufficiently large $C_0$, $C_5 > C_6$ holds and then, for all sufficiently large $C_0$,
\beql{eq-bound-diff2}
\E [g(X_{(C)}(s+)) - g(\hat X_{(C)}(s+))] I\{X_{(C)}(s+) \ge \hat X_{(C)}(s+)\} \le o(t).
\eeql
Observe also that, if $X_{(C)}(s+) \ge \hat X_{(C)}(s+)$, which is right after the single jump,
then 
\beql{eq-bound-diff222}
\mbox{$g(X_{(C)}(u)) - g(\hat X_{(C)}(u))$ can only decrease in the interval $[s+,t]$,}
\eeql
because there are no jumps in it.
From \eqn{eq-bound-diff2} and \eqn{eq-bound-diff222}, for all sufficiently large $C_0$,
\beql{eq-bound-diff4}
\E [g(X_{(C)}(t)) - g(\hat X_{(C)}(t))] I\{X_{(C)}(s+) \ge \hat X_{(C)}(s+)\} \le o(t).
\eeql
Further, if $\hat X_{(C)}(s+)>C_0$, and then necessarily $X_{(C)}(s+)=C_0$, we have
$$
g(X_{(C)}(s+)) < g(\hat X_{(C)}(s+))
$$
and then also
$$
g(X_{(C)}(t)) < g(\hat X_{(C)}(t)),
$$
because, recall, we only consider small $t$. Thus,
\beql{eq-bound-diff3}
[g(X_{(C)}(t)) - g(\hat X_{(C)}(t))] I\{X_{(C)}(s+) < \hat X_{(C)}(s+)\} \le 0 ~~\mbox{for small $t$}.
\eeql
Combining \eqn{eq-bound-diff4} and \eqn{eq-bound-diff3}, we finally obtain \eqn{eq-trans-bound}.
$\Box$

\subsection{Supplementary fact: Lemma~\ref{lem-conv-determ}}
\label{sec-suppl-fact}

The following supplementary fact, Lemma~\ref{lem-conv-determ}, 
 is used in the proof of 
Lemma~\ref{lem-gamma_B-tight}. It is used at the point where we assume that $C=C(B) \to \infty$ as $B\to\infty$, 
and the scaled distributions $\psi_{(C)}(y)=\gamma_B(Cy)$ are such that 
$\psi_{(C)} \stackrel{w}{\rightarrow} \psi_{(\infty)}$, where $\psi_{(\infty)}$ is a proper distribution;
it is also known that the speed $w_B$ converges, $w_B \to w_\infty$ (in fact, $w_\infty =v$).
Lemma~\ref{lem-conv-determ} basically says that as $C\to\infty$, the (fluid-scaled) process
$X_{(C)}(t)= (1/C) X_B(Ct)$, whose stationary distribution is $\psi_{(C)}$, in the limit becomes (informally speaking) a 
deterministic (fluid-limit) process $X_{(\infty)}(t)$, whose stationary distribution is $\psi_{(\infty)}$. 
The process $X_{(\infty)}(t)$
 is such that (informally speaking) 
$$
\frac{d}{dt} X_{(\infty)}(t) = \eta(\psi_{(\infty)}(X_{(\infty)}(t))) - w_\infty.
$$
As a result,
the distribution $\psi_{(\infty)}$ must be concentrated on those points $y$ where (informally speaking) 
$\eta(\psi_{(\infty)}(y))=w_\infty$.

\begin{lem}
%[Convergence of space/time-scaled processes to a deterministic process]
\label{lem-conv-determ}
Consider a subsequence of $B\to\infty$ and a function $C=C(B)$ such that $C\le B$, $C\to\infty$ and $B/C\to C_0^* \le \infty$. 
Assume that $w_B \to w_\infty$ for some constant $w_\infty \in (0,1)$. Denote $\psi_{(C)}(y) \doteq \gamma_B(Cy)$, and assume that $\psi_{(C)} \stackrel{w}{\rightarrow} \psi_{(\infty)}$, where $\psi_{(\infty)}$ is a proper distribution, not necessarily continuous. (It is concentrated on a finite interval $[-C_0^*,C_0^*]$ if $C_0^* < \infty$ or on the entire $\R$ otherwise.) For the chosen subsequence of $B$ and the corresponding values $C$, consider the sequence of processes $X_{(C)}(t)\doteq (1/C) X_B(Ct), t\ge 0,$ with converging initial states $X_{(C)}(0)\to y_0 \in \R$. 
(If $C_0^* < \infty$, then this necessarily means $y_0 \in [-C_0^*,C_0^*]$.)
Then, the following holds.

(i) The sequence of processes can be constructed on a common probability space so that, 
w.p.1, any further subsequence of $B$ has still further subsequence of $B$ along which $X_{(C)}(t) \to X_{(\infty)}(t)$, where
$X_{(\infty)}(t)$ is a deterministic Lipschitz continuous trajectory, such that $X_{(\infty)}(0)=y_0$ and for all $t\ge 0$,
\beql{eq-fluid-diff-incl1}
\frac{d^+_\ell}{dt} X_{(\infty)}(t) \ge \eta(\psi_{(\infty)}(X_{(\infty)}(t))) - w_\infty,~~ \mbox{if $X_{(\infty)}(t) < C_0^*$},
\eeql
\beql{eq-fluid-diff-incl2}
\frac{d^+_u}{dt} X_{(\infty)}(t) \le \eta(\psi_{(\infty)}(X_{(\infty)}(t)-)) - w_\infty, ~~ \mbox{if $X_{(\infty)}(t) > - C_0^*$}.
\eeql

(ii) Any limiting trajectory $X_{(\infty)}(\cdot)$ in (i) is such that
$X_{(\infty)}(t) \to [q^{low}, q^{up}]$, 
where 
$$
q^{low} = \sup\{y~|~\eta(\psi_{(\infty)}(y)) > w_\infty\}, ~~~ q^{up} = \inf\{y~|~\eta(\psi_{(\infty)}(y)) < w_\infty\}.
$$ 
Moreover, this convergence is uniform in $y_0$ and all such $X_{(\infty)}(\cdot)$, as long as 
 $y_0$ is restricted to a compact subset of $\R$.

(iii) Distribution $\psi_{(\infty)}$ is concentrated on the segment $[q^{low}, q^{up}]$.
Consequently, if $q^{low} < q^{up}$ (and then $\psi_{(\infty)}(y)=\eta^{-1}(w_\infty)$ for $y\in (q^{low},q^{up})$),
the distribution $\psi_{(\infty)}$ has exactly 
two atoms, at points $q^{low}$ and $q^{up}$, with masses 
$\eta^{-1}(w_\infty)$ and $1-\eta^{-1}(w_\infty)$, respectively. 
\end{lem}

We remark that Lemma~\ref{lem-conv-determ}(i) easily implies a stronger property: 
w.p.1 $X_{(C)}(\cdot) \stackrel{u.o.c.}{\rightarrow} X_{(\infty)}(\cdot)$, where $X_{(\infty)}$ 
is the unique Lipschitz trajectory, such that $X_{(\infty)}(0)=y_0$ and for almost all $t$ w.r.t. Lebesgue measure,
$$ 
X'_{(\infty)}(t) \in [\eta(\psi(X_{(\infty)}(t)))- w_\infty, \eta(\psi(X_{(\infty)}(t)-))- w_\infty].
$$ 
We do not need this stronger property in the present paper.

\begin{proof}[Proof of Lemma~\ref{lem-conv-determ}]
 The proof uses a fairly standard fluid-limit type argument. 
When $C$ is large, the process $X_{(C)}(\cdot)$ becomes ``almost deterministic.'' 
Indeed, when $C$ is large, and then the distribution $\psi_{(C)}$ is close to a proper distribution $\psi_{(\infty)}$,
the process makes right jumps at the rate $O(C)$ with the jump sizes being $O(1/C)$; as a result, 
in the vicinity of time $t$ the process trajectory 
is almost deterministic, with the derivative being close to $\eta(\psi_{(\infty)}(X_{(C)}(t))) - w_\infty$. All properties stated in the lemma
stem from this basic behavior.

We proceed with the formal proof.
Consider the following natural common probability space construction for the particle processes $X_B$, for all $B$. A unit rate 
Poisson process $\Pi(t), t\ge 0,$ drives the particle urges to jump. An i.i.d. sequence $Z_1, Z_2, \ldots,$ of random variables with distribution $J(\cdot)$ determines the sequence of particle jump sizes, when it does jump. An i.i.d. sequence $\Xi_1, \Xi_2, \ldots,$ of random variables,
uniformly distributed in $[0,1]$, determines whether or not the particle actually jumps when it gets the $i$-th urge to jump;
specifically, if a particle gets $i$-th urge to jump at time $t$, when its location is $y=X_B(t)$, it actually jumps
if $\Xi_i \le \bar \eta(y,\gamma_B)$. 
This construction ensures that, for each $B$, the process $X_B$ is (up to stochastic equivalence) as defined.

The driving sequences satisfy the functional strong law of large numbers (FSLLN) properties: as $B\to\infty$ (and 
corresponding $C=C(B) \to\infty$), w.p.1, 
\beql{eq-fslln1}
\left( \frac{1}{C} \Pi(Ct), ~t\ge 0 \right) \stackrel{u.o.c.}{\rightarrow} (t, ~t\ge 0), 
\eeql
\beql{eq-fslln2}
\left( \frac{1}{C} \sum_{i=1}^{\lfloor C s \rfloor} Z_i, ~s\ge 0 \right) \stackrel{u.o.c.}{\rightarrow} (s, ~s\ge 0), 
\eeql
\beql{eq-fslln3}
\left( \frac{1}{C} \sum_{i=1}^{\lfloor C s \rfloor} I\{\Xi_i \le \xi\}, ~0 \le \xi \le 1, ~s\ge 0 \right) 
\stackrel{u.o.c.}{\rightarrow} (\xi s, ~0 \le \xi \le 1, ~s\ge 0).
\eeql

(i) We have $X_{(C)}(t)=X_{(C)}(0) + X^{\uparrow}_{(C)}(t) - X^{\downarrow}_{(C)}(t)$, where 
$X^{\uparrow}_{(C)}(t)$ accounts for the forward (right) jumps in $[0,t]$ and $X^{\downarrow}_{(C)}(t)$ accounts for the total distance
travelled by the particle backwards (left) in $[0,t]$ (at the speed which may be either $w_B$ or $0$). Both $X^{\uparrow}_{(C)}$
and $X^{\downarrow}_{(C)}$ are non-decreasing; $X^{\downarrow}_{(C)}$ is $w_B$-Lipschitz; 
given the FSLLN properties \eqn{eq-fslln1}-\eqn{eq-fslln3}, we observe that the sequence of processes
$X^{\uparrow}_{(C)}$ 
is asymptotically Lipschitz, namely, w.p.1., for any $0\le t_1 \le t_2 < \infty$,
$$
\limsup_{C\to\infty} [X^{\uparrow}_{(C)}(t_2) - X^{\uparrow}_{(C)}(t_1)] \le t_2-t_1.
$$
This implies that, w.p.1., any subsequence of trajectories $X_{(C)}(\cdot)$ has a further subsequence along which
$$
X_{(C)}(\cdot) \stackrel{u.o.c.}{\rightarrow} X_{(\infty)}(\cdot),
$$
where $X_{(\infty)}(t)$ is $(1+w_\infty)$-Lipschitz, and $X_{(\infty)}(0)=y_0$.
It remains to show that a limit trajectory $X_{(\infty)}(\cdot)$ safisfies
\eqn{eq-fluid-diff-incl1}-\eqn{eq-fluid-diff-incl2}. This easily follows from \eqn{eq-fslln1}-\eqn{eq-fslln3}. 
Consider \eqn{eq-fluid-diff-incl1}, for example.
If $X_{(\infty)}(t)=y$, then for an arbitrarily small $\epsilon>0$ there exists $\delta>0$ such that
for all $t' \in (t-\delta,t+\delta)$ and all sufficiently large $C$ we have
$$
\psi_{(C)}(X_{(C)}(t')) < \psi_{(\infty)}(y)+\epsilon
$$
and then
$$
\bar \eta (X_{(C)}(t'), \psi_{(C)})  > \eta(\psi_{(\infty)}(y)+\epsilon).
$$
In other words, when $t'$ is close to $t$ and $C$ is sufficiently large, the probability of the particle of $X_{(C)}$ jumping when it gets an urge
is lower bounded by $\eta(\psi_{(\infty)}(X_{(\infty)}(t))+\epsilon)$. Application of \eqn{eq-fslln1}-\eqn{eq-fslln3} easily gives 
\eqn{eq-fluid-diff-incl1}. We omit further details. Property \eqn{eq-fluid-diff-incl2} is shown similarly.

(ii) Follows from (i).

(iii) Easily follows from (ii). Indeed, (ii) implies that for any compact set $K \in \R$ there exists $T>0$ such that, for all sufficiently large $C$, 
with uniformly high probability, if $X_{(C)}(0) \in K$ then $X_{(C)}(T)$ is close to $[q^{low}, q^{up}]$.
Therefore, the stationary distribution of $X_{(C)}(\cdot)$, i.e.  $\psi_{(\infty)}$, must in the limit concentrate on 
$[q^{low}, q^{up}]$.
\end{proof}

\subsection{Proof of Theorem~\ref{thm-existence}(ii): Finite moments} 

The distribution $\phi$ is proper and $\phi=A \phi$. Then \eqn{eq-moments} follows from Lemma~\ref{lem-stability},
\eqn{eq-moments-oper}. $\Box$

\subsection{Proof of Theorem~\ref{thm-main}} 

If assumption \eqn{eq-positive-density} holds, the mean-field model $f(x,t)$ is within the assumptions of that in \cite{GSS96}.
(Actually, in \cite{GSS96} it is assumed that $J(y) = 1-e^{-y}$ and $\eta(\nu)=(1-\nu)^K$, $K \ge 1$. However, all proofs in \cite{GSS96}
hold {\em as is} under more general assumption \eqn{eq-positive-density} on $J(\cdot)$, 
and for arbitrary continuous strictly decreasing $\eta(\cdot)$.) Theorem 1 of \cite{GSS96} states that {\em if a traveling wave shape $\phi$ with well-defined finite mean exists}, then statements (i) and (ii) of our Theorem~\ref{thm-main} hold. 
(Note that the uniform convergence in \eqn{eq-conv-to-shape} follows from the $L_1$-convergence, because $\phi(\cdot)$ is continuous non-decreasing and each $f(\cdot+vt,t)$ is non-decreasing.)
But, by 
our Theorem~\ref{thm-existence}, such traveling wave shape does exist, which completes the proof. 
$\Box$

\section{Discussion}
\label{sec-discusion}

\subsection{A generalization: more general jump process}

Our main results (along with proofs) easily extend, for example, to the following system. There are $n$ particles, moving 
in jumps in the positive direction. Each particle can make two types of jumps. Jumps of type $1$ are those in our original model:
they are driven by an independent Poisson process of rate $\mu \ge 0$,
 jump probabilities given by a non-increasing function $\eta_n(\cdot)$ uniformly converging to a
 strictly decreasing continuous function $\eta(\nu)$, and the jumps sizes are independent with
distribution $J(\cdot)$. Suppose that, in addition, each particle can make jumps of type $2$: they are driven by 
an independent Poisson process of rate $\mu_2 \ge 0$, at which points the particle jumps right w.p.1, and the jumps sizes are independent
with distribution $J_2(y), y\ge 0$; we denote $\bar J_2(y) = 1- J_2(y)$. Assume that \eqn{eq-moment-finite} holds for $J(\cdot)$, 
and the analogous assumption holds for $J_2(\cdot)$ as well, for the same integer $\ell \ge 2$ as in \eqn{eq-moment-finite}:
\beql{eq-moment-finite2}
m_2^{(k)} \doteq \int_0^\infty y^k dJ_2(y) < \infty, ~~k=1,2,\ldots,\ell.
\eeql
So, $m_2^{(1)}<\infty$ is the mean type $2$ jump size.

This generalization may be useful to model situations when each particle may change its state either independently (type $2$ jumps)
or depending on the locations of other particles (type $1$ jumps).

The corresponding mean-field model is still given by Definition~\ref{def-mfm}, except \eqn{eq-dyn-trans}
generalizes to
\beql{eq-dyn-trans-gen}
\frac{\partial}{\partial t} f(x,t) =
- \mu \int_{-\infty}^x d_y f(y,t) \bar \eta(y,f(\cdot,t)) \bar J(x-y) - \mu_2  \int_{-\infty}^x d_y f(y,t) \bar J_2(x-y).
\eeql
Speed $v$ generalizes to
$$
v \doteq \mu m^{(1)} \int_0^1 \eta(\xi) d\xi + \mu_2 m_2^{(1)},
$$
for which the ``conservation law'' \eqn{eq-speed-conserv} holds. Therefore, just as in the original model,
if the initial mean $\bar f(0)$ is well-defined finite, then it moves right at the constant speed $v$: $\bar f(t) = \bar f(0) + vt$.

A traveling wave shape is given by Definition~\ref{def-tws}, and easily shown to satisfy a more general form of \eqn{eq-wave},
\beql{eq-wave-gen}
v  \phi'(x) =
\mu \int_{-\infty}^x \phi'(y) \eta(\phi(y)) \bar J(x-y) dy  + \mu_2  \int_{-\infty}^x \phi'(y) \bar J_2(x-y) dy,
\eeql
for each $x$,  and the derivative $\phi'(x)$ is continuous.

For this more general mean-field model, {\em Theorem~\ref{thm-existence} (with condition \eqn{eq-moment-finite}
complemented by \eqn{eq-moment-finite2}) and Theorem~\ref{thm-main} 
(with condition \eqn{eq-positive-density} complemented by the analogous condition on $J_2(\cdot)$) 
hold as they are.} The proofs are same, up to straightforward adjustments. 

The key feature that the more general model shares with the original one -- and which makes the same analysis work -- 
is that  jumps of a particle depend only on its location quantile. (This, in particular, implies that
the speed of the mean $\bar f(t)$, and then of 
a traveling wave $\phi(x-vt)$, is known in advance and equal to $v$.) As long as this key feature is preserved for other mean-field models,
we believe that our main results and analysis have a good chance to extend to such other models as well.

\subsection{Condition \eqn{eq-positive-density} on the jump size distribution}

Condition \eqn{eq-positive-density} is used in  \cite{GSS96} to prove that the $L_1$-distance between any two
mean-field models (with equal mean) is strictly decreasing as long as these mean-field models are different.
(Actually, the results in \cite{GSS96} are specifically for the exponential $J(\cdot)$, but as far as the decreasing 
$L_1$-distance is concerned, only condition \eqn{eq-positive-density} on $J(\cdot)$ is used.)
The decreasing $L_1$-distance result of  \cite{GSS96} is then used in the proof of our Theorem~\ref{thm-main}.

It is likely that the analysis in \cite{GSS96} can be generalized to establish the decreasing $L_1$-distance property under a {\em much relaxed condition \eqn{eq-positive-density}.} If so, our Theorem~\ref{thm-main} holds under a relaxed condition \eqn{eq-positive-density}
as well.

\subsection{A conjecture about the limit of stationary distribution}

Consider the stochastic particle system, with $n$ particles. (Not the corresponding mean-field model, which is the focus of this paper.) 
Recall that by $f^n(\cdot,t)$ we denote its random state (the empirical distribution of particle locations) at time $t$. 
Denote by $f^n_*(\cdot,t)$ the function $f^n(\cdot,t)$, recentered so that its median is at $0$. 
Assume \eqn{eq-moment-finite} and \eqn{eq-positive-density}. It is not hard to see (using fluid limit technique, for example) that,
for any fixed $n$, the process $f^n_*(\cdot,t)$ is stochastically stable (positive Harris recurrent), and therefore has
unique stationary distribution. (We do not give details of the stability proof -- which are not hard -- because analysis of the stochastic particle
system with finite $n$ is not the subject of this paper.) By our Theorems~\ref{thm-existence} and \ref{thm-main}, there exists the unique
(up to a shift) traveling shape $\phi(\cdot)$.
 Then, the results of \cite{GMP97} (convergence 
of $f^n(\cdot,t)$ to a mean-field model $f(\cdot,t)$, as $n\to\infty$) and our Theorem~\ref{thm-main} 
(convergence of a mean-field model $f(\cdot,t)$ to the traveling wave, as $t\to\infty$) strongly suggest the following 
conjecture about the limit of stationary distributions of $f^n_*(\cdot,t)$. 

\begin{conj}
\label{conj1}
Assume \eqn{eq-moment-finite} and \eqn{eq-positive-density}.
Let $f^n_*(\cdot,\infty)$ be a random value of $f^n_*(\cdot,t)$ in the stationary regime.
Let $\phi(\cdot)$ be the unique traveling wave shape, centered to have the median at $0$. 
Then, as $n\to\infty$, $f^n_*(\cdot,\infty)$ concentrates at $\phi(\cdot)$, namely
$$
\|f^n_*(\cdot,\infty) - \phi(\cdot)\| \stackrel{P}{\rightarrow} 0.
$$
\end{conj}

As natural as it is, this conjecture (which is also supported by the simulation experiments in \cite{GSS96}) does {\em not} 
directly follow from the results of \cite{GMP97} and our Theorem~\ref{thm-main}. Establishing Conjecture~\ref{conj1} is a subject of future work.

\end{document}